\documentclass[10pt, reqno]{amsart}

\usepackage{amssymb,latexsym}
\usepackage{eucal}

\newtheorem*{Stt11}{Problem 1.1}
\newtheorem*{Stt12}{Theorem 1.2}
\newtheorem*{Stt13}{Theorem 1.3}
\newtheorem*{Stt14}{Theorem 1.4}
\newtheorem*{Stt15}{Theorem 1.5}
\newtheorem*{Stt16}{Theorem 1.6}
\newtheorem*{Stt17}{Theorem 1.7}
\newtheorem*{Stt19}{Theorem 1.8}
\newtheorem*{Stt110}{Corollary 1.9}
\newtheorem*{Stt111}{Theorem 1.10}
\newtheorem*{Stt112}{Corollary 1.11}
\newtheorem*{Stt113}{Problem 1.12}

\newtheorem*{Lem31}{Lemma 3.1}
\newtheorem*{Lem32}{Lemma 3.2}

\newtheorem*{Lem33}{Lemma 3.3}

\newcommand{\D}{\Delta}

\newcommand{\p}{\partial}
\newcommand{\e}{\varepsilon}

\newcommand{\A}{\mathcal{A}}
\newcommand{\ph}{\varphi}
\newcommand{\Ga}{\Gamma}

\newcommand{\R}{\mathcal{R}}
\newcommand{\PP}{\mathcal{P}}
\newcommand{\W}{\mathcal{W}}
\newcommand{\U}{\mathcal{U}}
\newcommand{\V}{\mathcal{V}}
\newcommand{\F}{\mathcal{F}}
\newcommand{\HH}{\mathcal{H}}
\newcommand{\N}{\mathcal{N}}
\newcommand{\s}{\mathcal{S}}

\newcommand{\ff}{{19/12}}

\newcommand{\E}{\text{E}}

\begin{document}

\title[On Dehn functions of infinite presentations of groups]
{On Dehn functions of infinite presentations of groups}
\author{R.\ I.\ Grigorchuk and S.\ V.\ Ivanov}

\address{Department of Mathematics\\   Mailstop 3368\\
     Texas A\&M University\\
     College Station\\ TX 77843-3368 \\ USA}
     \email{grigorch@math.tamu.edu}

\address{ Department of Mathematics\\
University of Illinois \\1409 West Green Street\\ Urbana\\   IL
61801\\ USA} \email{ivanov@math.uiuc.edu}

\thanks{The first named author is partially supported by NSF grants
DMS 04-56185, DMS 06-00975  and by the Swiss National Science Foundation.
The second named author is supported in part by NSF grant DMS 04-00746}

\begin{abstract}

We introduce two new types of  Dehn   functions of group presentations
which seem more  suitable  (than the  standard Dehn function)  for
infinite group presentations and prove the fundamental equivalence
between the solvability of the word problem for a  group presentation
defined by a decidable set of defining words and the property of being
computable for one of the newly introduced functions (this equivalence
fails for the standard Dehn function). Elaborating on this equivalence
and making use of this function, we obtain a characterization of finitely
generated groups for which the word problem can  be solved in
nondeterministic polynomial time.

We also give upper bounds for these functions, as well as for the
standard Dehn function, for two well-known periodic groups. In
particular, we prove that the (standard) Dehn function of a 2-group $\Ga$
of intermediate growth, defined by a system of defining relators due to
Lysenok, is bounded from above by $C_1 x^2\log_2x$, where $C_1>1$ is a
constant. We also show that the (standard) Dehn function of a free
$m$-generator Burnside group $B(m,n)$ of exponent $n \ge 2^{48}$, where
$n$ is either odd or divisible by $2^9$, defined by a minimal system of
defining relators, is bounded from above by the subquadratic function
$x^{\ff}$.

\end{abstract}
\maketitle

\section{Introduction}

Let a finitely generated group $G$ be defined by a presentation
in terms of generators and defining relators
\begin{equation}\label{pr3}
G  = \langle \ \A \ \| \ \R  \ \rangle \ ,
\end{equation}
where $\A = \{ a_1, \dots , a_m\}$ is a finite alphabet  and $\R$ is a
set of defining { relators} which are nonempty cyclically reduced words
over the alphabet $\A^{\pm 1} = \A \cup \A^{-1}$. Let $F(\A)$ denote the
free group over $\A$, $|W|$ mean  the length of a word $W \in F(\A)$ over
the alphabet $\A^{\pm 1}$, and $\langle\langle \R \rangle\rangle $ denote
the normal closure of $\R$ in $F(\A)$. Then the notation \eqref{pr3}
means that $G$ is the quotient group $F(\A)/\langle\langle \R
\rangle\rangle $. Recall that a presentation \eqref{pr3} is called {\em
finite} if  $\R$ is finite and $\R$ is termed {\em decidable} (or
recursive)  if there is an algorithm to decide whether a given word over
$\A^{\pm 1}$ belongs to $\R$.

Let $K(\A, \R)$ be a  2-complex  associated with the presentation
\eqref{pr3} so that $K(\A, \R)$ has a single 0-cell, oriented 1-cells of
$K(\A, \R)$ are in bijective correspondence with letters of $\A^{\pm 1}$,
and 2-cells of $K(\A, \R)$ are in bijective correspondence with the words
of $\R$ that  naturally determine the attaching maps of the 2-cells. Thus
$K(\A, \R)$ is the standard geometric realization of \eqref{pr3} with the
fundamental group $\pi_1(K(\A, \R))$ being isomorphic to~$G$.

By a {\em van Kampen diagram} over the presentation \eqref{pr3} we mean a
planar, finite, connected and simply connected 2-complex $\D$ which is
equipped  with a continuous cellular map $\mu : \D \mapsto K(\A, \R)$
whose restriction on every cell of $\D$ is a homeomorphism. By an {\em
edge} of $\D$ we mean the closure of a 1-cell.  If $e$ is an oriented
edge  of $\D$ and $\mu(e)$ corresponds to a letter $a \in \A^{\pm 1}$,
then $a$ is termed the {\em label} of $e$ and is denoted $\ph(e)$. Note
$\ph(e^{-1})= \ph(e)^{-1}$, where $e^{-1}$ denotes the edge with opposite
orientation. If $p =e_1 \dots e_k$ is a path in $\D$, where $e_1, \dots,
e_k$ are oriented edges of $\D$, then we set $\ph(p) {}= \ph(e_1) \dots
\ph(e_k)$. According to a well-known lemma of van Kampen, see \cite{LS,
Ol89}, a word $W$ belongs to $\langle\langle \R \rangle\rangle $ if and
only if there exists a van Kampen diagram over \eqref{pr3} whose boundary
path $\p \D$ is labeled by the word $W$, in which case we write $\ph(\p
\D ) \equiv W$, where the sign $\equiv$ means the literal equality of
(cyclic) words. The number of $j$-cells of a van Kampen diagram $\D$ is
denoted by $| \D(j)|$, $j \in \{0,1,2\}$.

Let $W \in \langle\langle \R \rangle\rangle$ and $j \in \{
0,1,2\}$.  Define $L_j(W)$ to be the minimal number of
$j$-cells in a van Kampen diagram $\D$ over \eqref{pr3} whose
boundary $\p \D$ is labeled by the cyclic word $W$, that is
\begin{equation} \label{Lj}
L_j(W) {}= \min \{ \,  | \D(j) | \mid \ph(\p \D ) \equiv W \} \
.
\end{equation}
For an integer $x \ge 1$, define
\begin{equation} \label{DFj}
f_j(x) {}= \max \{ \,  L_j(W) \mid W \in  \langle\langle \R
\rangle\rangle \  \text{and} \ |W| \le x  \} \ .
\end{equation}
 Recall that if $W \in
\langle\langle \R \rangle\rangle $, then, in $F(\A)$, one has
an equality of the form
\begin{equation}\label{FD}
W = \prod_{i=1}^L X_i R_i^{\e_i} X_i^{-1} \ ,
\end{equation}
where $X_i \in F(\A)$, $R_i \in \R$, $\e_i \in \{\pm 1\}$ and
where we allow $L =0$ (when $W = 1$ in  $F(\A)$). An equivalent
way to define the number $L_2(W)$ is to pick minimal $L \ge 0$
over all products of the form \eqref{FD}.

When $j=2$,  the foregoing definition \eqref{DFj} defines  the
well-known {\em  Dehn function} $f_2(x)$ of a group
presentation \eqref{pr3} that  has been subject of intensive
research for the past twenty years. Recall that the concept of
the Dehn function $f_2(x)$ of a group presentation \eqref{pr3}
was introduced into group theory by Gromov in his seminal
article \cite{Gromov} in 1987 and now is fundamental in
geometric group theory.  For instance, a finite presentation
\eqref{pr3} defines a (word) hyperbolic group if and only if
its Dehn function $f_2(x)$ is bounded from above by a linear
function. For a finite presentation \eqref{pr3}, the
solvability of the word problem is equivalent to the property
of being computable for the function $f_2(x)$. As was pointed
out by the referee, in 1985, Madlener and Otto \cite{MaOt}
considered a notion of the ``derivational complexity" which was
an earlier version of the definition of the Dehn function and
which was defined by means of regarding a presentation of a
monoid, or a group, as a string rewriting system.

The object of this paper is to introduce and study the functions
$f_0(x)$, $f_1(x)$, to establish basic relations between the solvability
of the word problem for presentation \eqref{pr3} and the property of
being computable for functions $f_0(x)$, $f_1(x)$,  and to give upper
bounds for the functions $f_j(x)$ in case of  presentations of two
well-known periodic (or torsion) groups. The function $f_j(x)$, $j \in
\{0,1,2\}$, will be referred to as the {\em Dehn $j$-function} of a group
presentation \eqref{pr3} (with the prefix ``$j$-'' frequently omitted).

As was pointed out above, for  a finite presentation
\eqref{pr3}, the solvability of the word problem is equivalent
to the property of being computable for the Dehn function
$f_2(x)$. However, if $\R$ is not finite, then this equivalence
breaks down and generates  a few interesting  problems which we
state for two  Dehn functions $f_2(x)$ and $f_0(x)$.

\begin{Stt11} Let the relator set $\R$ of
a presentation \eqref{pr3} be decidable and $j \in \{ 0, 2\}$.
Prove or disprove that

$(\mathrm{a})$  If the word problem for \eqref{pr3} is
solvable,
 then the Dehn $j$-function $f_j(x)$  of \eqref{pr3}  is
 computable.

$(\mathrm{b})$   If the Dehn $j$-function  $f_j(x)$  of
\eqref{pr3} is computable,  then the word problem for
\eqref{pr3} is solvable.
\end{Stt11}

Remarkably,  the Dehn $1$-function $f_1(x)$ can be used in place of
$f_2(x)$ to fix the  failing equivalence, see Example~2.4. Taking
advantage of counting 1-cells in place of 2-cells, we have the following
basic result.

\begin{Stt12} Let  $\R$ in \eqref{pr3} be  decidable.
Then the word problem for \eqref{pr3} is solvable if and only
if the Dehn $1$-function $f_1(x)$ of \eqref{pr3} is computable.
\end{Stt12}

Elaborating on this equivalence, we obtain a characterization
of finitely generated groups for which the word problem could
be solved in nondeterministic polynomial time.

\begin{Stt13} Let a group $G = \langle a_1, \dots, a_m \rangle$
be generated by elements $a_1, \dots, a_m$. Then the word
problem for $G$ is in {\bf NP}, i.e. it can be solved
algorithmically in nondeterministic polynomial time,  if and
only if there exists a presentation $\langle a_1, \dots, a_m \,
\| \,  \R \rangle$ for $G$  such that its  Dehn $1$-function
$f_1(x)$ is bounded by a polynomial and the problem to decide
whether a word $W = W(a_1^{\pm 1}, \dots, a_m^{\pm 1})$ belongs
to $\R$ is in {\bf NP}.
\end{Stt13}

Let us emphasize that the complexity, as well as the solvability, of the
word problem for a finitely generated group $G$, given by a presentation
\eqref{pr3}, does not actually depend on $\R$. Indeed, we need to decide
whether or not a given word $W$ over $\A^{\pm 1}$  represents the
identity element of $G$, i.e. whether  $W \in \langle \langle \R \rangle
\rangle$. On the other hand, when speaking about the word problem for a
finitely generated group $G$, we always assume that a finite generating
set, say $\A$, is fixed for $G$, i.e. $G$ is regarded as the quotient
group $F(\A)/\langle \langle \R \rangle \rangle$ for some $\R$, see, for
example, Theorem~1.3 above and Corollary~1.11 below.
\medskip

Theorem~1.3 is reminiscent of a deep result of Birget, Ol'shanskii, Rips
and Sapir~\cite{BORS} that states that a finitely generated group $G$ has
the word problem in {\bf NP} if and only if $G$ is isomorphic to a
subgroup of a finitely presented group whose Dehn function is bounded by
a polynomial. However, unlike the result of \cite{BORS}, Theorem~1.3 is a
straightforward  corollary of the advantageous definition of the
 function $f_1(x)$. Unlike the result of
\cite{BORS}, Theorem~1.3 also holds for other computational classes, for
example, it holds with {\bf PSPACE}   in place of {\bf NP} (recall that
the class {\bf PSPACE} consists of decision problems that could be solved
in polynomial space, for more details see \cite{Pap}).
\medskip

Let $\F$ denote the set of functions $f : \mathbb N \to \mathbb N \cup \{
0 \}$, where $\mathbb N = \{ 1,2, \dots \}$ is the set of natural
numbers. If $f, g \in \F$, we write $f \preccurlyeq  g$ if there is an
integer $C > 0$ such that $f(x) \le C g(C x) +Cx$ for every $x \in
\mathbb N$. Two functions $f$, $g \in \F$ are said to be (linearly)
equivalent, denoted $f \simeq g$, if $f \preccurlyeq  g$ and  $g
\preccurlyeq f$. It is evident that this relation $\simeq$ is indeed an
equivalence relation and the equivalence class $[ f ]_\simeq$ could be
regarded as the growth rate of a function $f \in \F$.

The natural questions on relations between functions $f_j(x)$,
$j=0,1,2$,  are addressed in the following.

\begin{Stt14}
$(\mathrm{a})$   For every presentation~\eqref{pr3},
 $ f_0(x)  \le 2f_1(x)$ and  $f_2(x)  \le 2f_1(x)$. In particular,
 $ f_0(x)  \preccurlyeq f_1(x)$ and $f_2(x)  \preccurlyeq
 f_1(x)$.

$(\mathrm{b})$ Let every element in $\R$ of \eqref{pr3} be a
word of length $> 1$. Then the Dehn functions $f_0(x)$ and
$f_1(x)$ of \eqref{pr3} are equivalent.
\end{Stt14}

Making use of the ideas of proofs of Theorems~1.2 and 1.4(b),
we will give a  positive solution to Problem~1.1 in a special
case for $j =0$.

\begin{Stt15}
 Let $\R$ in \eqref{pr3} be decidable and, for
every  $R \in \R$,   $|R| > 1$. Then the word problem for
\eqref{pr3} is solvable if and only if the Dehn $0$-function
$f_0(x)$ of \eqref{pr3} is computable. Moreover,  $f_0(x)$ is
computable if and only if $f_1(x)$ is computable.
\end{Stt15}

\medskip

To generalize the Andrews--Curtis and Magnus conjectures, the
following operation over  a finite presentation \eqref{pr3} is
introduced in the  article \cite{Iv06}: An element of $\R$ is
replaced by  another element if doing so does not change the
normal closure of $\R$. More generally, consider replacement of
a finite subset $\s$ of $\R$  by  another finite set $\U$ of
cyclically reduced words if doing so does not change the normal
closure of $\R$. As in \cite{Iv06}, call this operation a {\em
$T$-transformation}. Another operation over a presentation
\eqref{pr3}, called {\em stabilization}, is to add/delete a
letter $b$   to/from both $\A$ and $\R$ (when deleting $b$ from
$\A$ and $\R$, $b^{\pm 1}$  must not occur in any other word of
$\R$). Note that both $T$-transformation and stabilization are
special cases of Tietze transformations, see \cite[Section
1.5]{MKS}.

\begin{Stt16}  $(\mathrm{a})$  \ Let $\langle \A \,  \| \, \R'\rangle$
be obtained from a presentation \eqref{pr3}  by a
$T$-transformation and $f'_j(x)$, $f_j(x)$, $j=0,1,2$,  be
their corresponding Dehn functions. Then $\quad f'_0(x) \preccurlyeq
f_1(x)$, \ $f'_1(x) \simeq f_1(x)$, \ $f'_2(x) \simeq f_2(x)$.

$(\mathrm{b})$ \ Let  $\langle \A' \,  \| \, \R'\rangle$ be
obtained from a presentation \eqref{pr3} by a stabilization and
$f'_j(x)$, $f_j(x)$, $j=0,1,2$,  be their corresponding Dehn
functions. Then $f'_j(x) \simeq f_j(x)$ for $j =0,1,2$.
\end{Stt16}

Making use of $T$-transformations and stabilizations, we will
prove the equivalence  between all of the Dehn $j$-functions
$f_j(x)$, $j=0,1,2$, for a finite presentation.

\begin{Stt17}
 Let \eqref{pr3} be a finite presentation. Then all of
its Dehn  functions $f_0(x)$, $f_1(x)$, $f_2(x)$ are
equivalent.
\end{Stt17}

In Section 2, we will give Examples~2.1--2.2 of group
presentations for which the pairs $(f_1(x), f_2(x))$, $(f_0(x),
f_1(x))$, $( f_0(x), f_2(x))$ consist of  nonequivalent
functions. Hence, in general, $f_1(x) \not\simeq f_2(x)$,
$f_0(x) \not\simeq f_1(x)$,  $f_0(x) \not\simeq f_2(x)$. In
Section 2, we will also consider Examples~2.3--2.4 of  group
presentations that address parts (a)--(b) of Problem~1.1 in the
case $j=2$.
\medskip

It is well known (e.g. see \cite{BMS}, \cite{MaOt}) that if $\PP$, $\PP'$
are two finite presentations of a group $G$ and $f_2(x)$, $f'_2(x)$ are
their corresponding Dehn $2$-functions, then $f_2(x) \simeq  f'_2(x)$
(this also follows from Theorem~1.6(a)--(b) because $\PP'$ can be
obtained from $\PP$ by stabilizations and $T$-transformations, see
\cite[Section 1.5]{MKS}). Furthermore, if $H$ is a subgroup of finite
index of a group $G$ that has a finite presentation (or, more generally,
$G$ and $H$ are abstractly commensurable), then the Dehn $2$-functions of
finite presentations of $G$ and $H$ are equivalent. We also recall that,
if the Dehn $2$-function of a finite presentation of a group $G$ is
subquadratic, then the group $G$ is word hyperbolic, see \cite{Gromov,
Ol91}. For more results on Dehn $2$-functions see \cite{Br, BB, BORS,
Guba, GS, OS2, OS1, OS4} and the references cited there.

On the other hand, if we allow infinite relator sets $\R$ in
\eqref{pr3}, then the growth rates $[f_j(x)]_\simeq$, $j
=0,1,2$, of  the Dehn $j$-function  of $G$ are  no longer
independent of a presentation of $G$. Indeed, if $\R$ consists
of all cyclically reduced words of $\langle \langle \R \rangle
\rangle$, then $f_2(x) \le 1$, $f_0(x) \le x$, $f_1(x) \le x$,
and hence $f_j(x) \simeq  x$, $j =0,1,2$.

In general, an inclusion $\R' \subseteq \R$, provided  $\langle \langle
\R'  \rangle \rangle =  \langle \langle \R  \rangle \rangle$,  easily
implies the inequalities $f_j(x) \le f_j'(x)$, $j =0,1,2$, for the
corresponding    Dehn $j$-functions (cf.  Theorem~1.6(a)) and, for this
reason, it is more natural to consider presentations \eqref{pr3} with
{\em minimal} $\R$, i.e., if  $\R' \subseteq \R$ and   $\langle \langle
\R'  \rangle \rangle =  \langle \langle \R  \rangle \rangle$ then $\R' =
\R$.  However, Y. de Conrulier (private communication) pointed out to us
that there are  finitely generated groups that possess no  presentations
\eqref{pr3} with minimal $\R$.

Investigation of  Dehn $j$-functions $f_j(x)$, $j =0,1,2$, of a
presentation \eqref{pr3}, where the relator set $\R$  need not be finite,
seems to be an interesting and important problem, especially in the case
when $\R$ is minimal. In this article, we obtain two results in this
direction for infinite presentations of periodic groups investigated by
the authors in earlier articles.

\medskip

Let $\Ga$ denote the 2-group of intermediate growth that was
originally discovered by the first author in \cite{G80} and
later investigated in \cite{G84}, \cite{Lys85} and other
papers, see also  \cite[ Chapter VIII]{harpe}.
It was shown by Lysenok \cite{Lys85}
that the group $\Ga$ can be defined by the following
presentation in terms of generators and defining relators
\begin{align}\notag
\Ga   = \langle \ a,b,c,d \ \| & \ a^2, \ b^2,\ c^2,\ d^2, \
bcd, \  {\sigma}^i((ad)^4), \\ & {\sigma}^i((adacac)^4), \ i
\geq 0 \ \rangle \ = \  \langle \ a,b,c,d{{}} \  \| \
\R(\infty) \ \rangle \ , \label{pr1}
\end{align}
where it is assumed that $\sigma^0 = \mbox{id}$ is the identity
map and that $\sigma$ is the endomorphism of the free group
$F(a,b,c,d)$ with the basis $\{ a,b,c,d\}$ defined by
\begin{align}\label{sig}
\sigma {}=  \begin{cases} a
 \mapsto aca \ , \\  b \mapsto  d \ , \\
c \mapsto b \ ,
\\ d \mapsto  c \ .
\end{cases}
\end{align}
The  group  $\Ga$  has no  finite  presentation, see \cite{G84, G99},
and, as is shown in \cite{G99}, the  set $\R(\infty)$ of  Lysenok's
relators is {minimal}.

\begin{Stt19}\label{th11} The Dehn $j$-functions
$f_{j, \Ga}(x)$,  \ $j = 0,1,2$,  of the Lysenok presentation
\eqref{pr1} of the $2$-group $\Ga$ are bounded from above by
$C_{1} x^2\log_2x$, where $C_{1} > 1$ is a constant.
\end{Stt19}

It is of interest to point out that, when proving Theorem~1.8, we will
construct a larger relator set $\R^*(\infty)$, containing $\R(\infty)$
(so $\R^*(\infty)$ is not minimal), and obtain different upper bounds for
the functions $f_{j, \Ga}^*(x)$, $j \in \{ 1,2\}$, corresponding to
$\R^*(\infty)$: $f_{1, {\Ga}}^*(x) < C^*_1 x^2 \log_2x$ and $f_{2,
{\Ga}}^*(x) < C^*_1x^2$, where $C^*_1 > 1$ is a constant.

As in article \cite{G98}, consider  a group defined by the
following presentation
\begin{align*}\notag
\Ga_t   = \langle \  a, b, c, d,   t  \ \| \  \R(\infty), \ \ t g t^{-1}
\sigma(g)^{-1} , \ g \in \{ a,b,c,d \}  \  \rangle \ .
\end{align*}
Since $\sigma(R) \in \langle \langle \R(\infty) \rangle
\rangle$ for every $R \in \R(\infty)$, it follows that $\Ga_t $
is an ascending HNN-extension of $\Ga$ with the stable letter
$t$. Thanks to the form of Lysenok's relators $\R(\infty)$, the
group $\Ga_t$ can also be defined by the following finite
presentation
\begin{align}\notag
\Ga_t   = \langle   \ a,b,c,d, t \ \| \ a^2, \ b^2, & \ c^2,\
d^2,  \ bcd, \  (ad)^4,  \ (adacac)^4 , \  \\ &
 \ t g t^{-1} \sigma(g)^{-1} , \ g \in \{ a,b,c,d \}     \  \rangle \
 .
 \label{pr1a}
\end{align}
It is immediate from the definition that $\Ga_t $ is a finitely presented
torsion-by-cyclic group, see also \cite{Iv05}, \cite{OS} for more
examples of  finitely presented torsion-by-cyclic groups. Furthermore, as
was observed in \cite{G98},  the group $\Ga_t$ is amenable but not
elementary amenable, a property shared by $\Ga$. As a consequence of
Theorem~1.8, we will obtain

\begin{Stt110}\label{cor1} The Dehn   $j$-functions
$f_{j, \Ga_t}(x)$,  \ $j = 0,1,2$, of
the finite presentation \eqref{pr1a} are  bounded from above by
$C_{2} 2^x x$, where $C_{2}>1$ is a constant.
\end{Stt110}

\medskip

Recall that a free $m$-generator Burnside group $B(m,n)$ of
exponent $n$ is the quotient $F/F^n$, where $F$ is a free group
of rank $m$ and $F^n = \langle W^n \mid W \in F \rangle$.

To construct a presentation for $B(m,n)$, as in \cite{ Iv94,
Ol82}, we consider a total order $\preceq$ on the set of all
words over the alphabet $\A^{\pm 1} = \{a_1^{\pm 1}, \dots,
a_m^{\pm 1}\}$ such that $ U \preceq V$ when $|U| \le |V|$. As
above, let $F(\A)$ denote the free group with the basis $\A$.
Set $B(m,n, 0) {}= F(\A) = \langle \A \ \| \ \emptyset
\rangle$. Proceeding by induction on $i \ge 1$, assume that the
group presentation $B(m,n, i-1)$ is already constructed. Let
$A_i$ be the minimal (if it exits), relative to the order
$\preceq$, word over $\A^{\pm 1}$  such that the image of $A_i$
has infinite order in the group defined by $B(m,n, i-1)$. Note
that $A_i$ may not exist and then our inductive process
terminates and results in the presentation $B(m,n, i-1)$. If
$A_i$ does exist, then the presentation $B(m,n, i)$  is
obtained from $B(m,n, i-1)$ by addition of the relator $A_i^n$.
Clearly,
\begin{equation}\label{pr2i}
B(m,n,i)   = \langle \ \A  \ \| \ A_1^n, \  A_2^n, \  \dots, \
A_i^n \  \rangle
\end{equation}
and, if $A_i$  exists for every $i \ge 1$, then, taking the
limit, we obtain
\begin{equation}\label{pr2}
B(m,n,\infty)   = \langle \ \A  \ \| \ A_1^n, \  A_2^n, \
\dots, \  A_i^n  , \ \dots \    \rangle \ .
\end{equation}

Assume that $n \ge 2^{48}$ and either $n$ is odd or divisible by $2^9$.
Under this assumption, it is proved in \cite[Theorem B]{Iv94}, see also
\cite{Iv92}, \cite{Iv98}, that $A_i$ does exist\footnote[1]{In the case
when $B(m,n)$ is finite, the word $A_i$ will fail to exist for some $i$.
We do not know whether it is possible that $B(m,n)$ is infinite and $A_i$
does not exist for some $i$.} for every $i \ge 1$ and the limit group
$B(m,n,\infty)$ is naturally isomorphic to the free Burnside group
$B(m,n)= F(\A)/F(\A)^n $ with the basis $\A$. Furthermore, it is shown in
\cite{Iv94}, see Theorem B and Lemma 21.1,  that the relator set $\{
A_1^n, \ A_2^n, \dots, A_i^n, \dots\}$ is decidable and minimal. We
remark that, for odd $n
> 10^{10}$, these results are due to Ol'shanskii \cite{Ol82,
Ol89}, compare with Novikov--Adian's \cite{NA}, Adian's
\cite{A75}, and Lysenok's \cite{Lys96}  presentations and
results on $B(m,n)$.

\medskip

For free Burnside groups, we  prove

\begin{Stt111} \label{Stt110} Let $m \ge 2$, $n \ge 2^{48}$ and
 $n$ be either  odd or divisible by $2^9$. Then the Dehn
$1$-function $f_{1, B}(x)$ of the presentation \eqref{pr2}  of a free
$m$-generator Burnside group $B(m,n)$ of exponent $n$ is bounded from
above by the subquadratic function $x^{\ff}$. In addition, $f_{0, B}(x)
\le   2 x^{\ff}$ and $f_{2, B}(x) \le \frac 2 n x^{\ff}$. Furthermore,
for every $i \ge 0$, the same upper bounds hold for the Dehn
$j$-functions  $f_{j, B(i)}(x)$ of the finite presentation $B(m,n,i)$
defined by \eqref{pr2i}.
\end{Stt111}

Note  that the cubic upper bound $f_{2, B}(x) \le 6(n^{-1}x)^3$ for the
Dehn function $f_{2, B}(x)$ of $B(m,n, \infty)$, in the case of odd $n >
10^{10}$, is due to Storozhev \cite[Section 28.2]{Ol89}. We also mention
that a linear  bound $f_{2, B(i)}(x) \le K_i x$ for the Dehn function
$f_{2, B(i)}(x)$  of the finite presentation $B(m,n,i)$ was obtained in
\cite[Lemma~21.1]{Iv94}. This linear upper bound implies that, for every
$i \ge 0$, the group $B(m,n,i)$ is word  hyperbolic. However, the
constant $K_i$, as a function of $i$, grows exponentially and so the
bounds $f_{2, B(i)}(x) \le K_i x$, $i =1,2,\dots$, do not shed any light
upon an upper bound  for $f_{2,B}(x)$. On the other hand, in view of
results of Gromov \cite{Gromov} and Ol'shanskii \cite{Ol91} on the
hyperbolicity of finitely presented groups with subquadratic Dehn
2-function, the subquadratic bound $f_{2, B(i)}(x) \le \frac 2 n x^{\ff}$
of Theorem~1.10 implies that, for every $i \ge 0$, the group
 given by presentation~\eqref{pr2i}  is word
hyperbolic and hence, by Theorem~1.7, its Dehn $j$-functions $f_{j,
B(i)}$ are bounded by a linear function.

As a consequence of Theorem~1.10 and  lemmas of \cite{Iv94}, we will also
derive

\begin{Stt112} Let $n \ge 2^{48}$ and either $n$ be odd or
divisible by $2^9$. Then the word and conjugacy problems for
the free Burnside group $B(m,n) = F(\A)/F(\A)^n$ are in $\bf
{NP}$.
\end{Stt112}

For odd $n > 10^{10}$, Corollary~1.11 could be derived from the
Storozhev's cubic bound  for $f_{2,B}(x)$ and lemmas of \cite{Ol82}. It
would be interesting to further investigate the complexity of the word
problem for $B(m,n)$ and find out whether it is in $\bf {P}$, i.e. it is
solvable in deterministic polynomial time, or in $\bf {coNP}$ or,
perhaps, $\bf {NP}$-complete, see \cite{Pap}. We also remark that the
word problem for the group $\Ga$, see \eqref{pr1}, is known to be in {\bf
P}. In fact, the word problem for $\Ga$ can be solved in subquadratic
time ($\sim y\log_2y$)), see \cite{G99}. A deterministic algorithm for
the conjugacy problem for
 $\Ga$, whose running time is at least exponential, as
well as the  history  of  the  question and further references
could be found in \cite{G05}. It would also be desirable to make
further progress on the following.

\begin{Stt113} Obtain  nontrivial lower bounds
for the Dehn $j$-functions $f_{j,\Ga}$, $f_{j,B}$, $j \in
\{1,2\}$, of presentations \eqref{pr1}, \eqref{pr2}, improve on
upper bounds for   $f_{j,\Ga}$, $f_{j,B}$, and determine the
growth rates $[f_{j,\Ga}]_\simeq$, \ $[f_{j,B}]_\simeq$.
\end{Stt113}

Note that, by Theorem~{1.4}(b), $f_{0,\Ga} \simeq f_{1,\Ga}$
and  $f_{0,B} \simeq f_{1,B}$.
\medskip

An interesting notion of the verbal Dehn function $f_w(x)$ of a variety
of groups, defined by a single identity $w =  1$, was introduced and
investigated by Ol'shanskii and Sapir \cite{OS2}. To give the definition,
we let $F_\infty=F(a_1, a_2, \dots)$ be the free group over the countably
infinite alphabet $\A_\infty = \{ a_1, a_2, \dots \}$, $w$ be a word, and
$w(F_\infty)$ be the $w$-verbal subgroup of $F_\infty$ generated by all
values of $w$ on $F_\infty$.  For every $U \in w(F_\infty)$, consider a
product in  $F_\infty$  of the form
 $$
 U = \prod_{j=1}^L
X_j w(Y_{1j}, \dots, Y_{kj})^{\e_j} X_j^{-1} \ ,
$$
where $X_j, Y_{i j} \in  F_\infty$  and $\e_j = \pm 1$. Taking
the minimal sum $\sum_{j=1}^L (|Y_{1j}| + \dots +|Y_{kj}|)$
over all such products for $U$, we obtain a number $s(U)$. Then
the {\em  verbal Dehn function} for the word $w$ is defined by
 $f_w(x) {}= \max \{ s(U) \mid |U| \le x, \ U \in w(F_\infty) \}$.

According to \cite{OS2}, in an unpublished work, for odd $n> 10^{10}$,
Mikhailov gave an upper bound  $x^{1+\e_n}$, where $\e_n >0$ and $\e_n
\to 0$ as $n \to \infty$,   for the verbal Dehn function $f_{z^n}(x)$ of
the Burnside variety of groups of exponent $n$.  We conjecture that an
analogous upper  bound holds for the function $f_{1,B}(x)$, where $n \gg
1$ is defined as in Theorem~1.10. However, it remains unclear whether
there are any relations between the functions $f_{ z^n}(x)$ and
$f_{1,B}(x)$ and whether an upper  bound   for $f_{z^n}(x)$ could
possibly yield any bound for $f_{1,B}(x)$.

As is pointed out in \cite{OS2}, the verbal Dehn function
$f_w(x)$ is {\em superadditive} for every word $w$, i.e.
$f_w(x)$ satisfies the inequality $f_w(x_1 +x_2) \ge f_w(x_1) +
f_w(x_2)$. Recall that the superadditive closure $\bar f(x)$ of
a function $f(x)$, $x \in \mathbb N$, is $\bar f(x) {}= \max \{
f(x_1) + \dots + f(x_r) \}$ over all sums $x =x_1 + \dots +
x_r$. It is immediate that  $\bar f(x)$ is superadditive.  A
conjecture, put forward by Guba and Sapir \cite{GS}, claims
that the Dehn function $f_2(x)$ of a finite presentation is
equivalent to its superadditive closure $\bar f_2(x)$. It seems
to be of interest to ask the same question for Dehn
$j$-functions $f_j(x)$,  $j=0,1,2$, defined for an arbitrary
presentation \eqref{pr3}, that is, to ask whether $f_j(x)
\simeq \bar f_j(x)$, $j=0,1,2$. Speaking of the  equivalence
$f_j(x) \simeq \bar f_j(x)$, we remark that we do not know
whether functions $f_{j,B}(x)$, $f_{j, \Ga}(x)$ are equivalent
to their superadditive closures. We do not know either  whether
two functions $f_{1,B}(x)$, corresponding to two different
choices of defining words $A_1^n$, $A_2^n, \dots$ in
\eqref{pr2}, are equivalent, to leave alone the equivalence of
functions $f_{1,B}(x)$ of  free Burnside groups $B(m,n)$ of
different ranks $m \ge 2$ and exponents $n \gg 1$.

\section{Proofs of Theorems~1.2--1.7}

{\em Proof of Theorem~1.2.}  Suppose that the word problem is
solvable for a  presentation \eqref{pr3} with  a decidable set
$\R$. Then the normal closure $\langle \langle  \R \rangle
\rangle$ of $\R$ is a decidable subset of $F(\A)$. To compute
the value of $f_1(x)$ for a given integer $x \ge 1$, consider
the set $\W_x$ of all words over the alphabet $\A^{\pm 1}$ of
length $\le x$. For each $U \in \W_x$, we determine whether $U
\in \langle \langle  \R \rangle \rangle$. If so, we construct a
van Kampen diagram $\D_U$  over \eqref{pr3} with $ \ph(\p \D_U)
\equiv U$. Now we check all van Kampen diagrams $\D$ such that
$\ph(\p \D) \equiv U$ and $|\D(1)| \le | \D_U(1)|$  and find a
diagram $\D_{U, 1}$ over \eqref{pr3}  with the minimal number
$|\D_{U,1}(1)|$. Then $f_1(x) = \max \{ |\D_{U,1}(1)| = L_1(U)
\mid U \in   \W_x \cap  \langle \langle  \R \rangle \rangle
 \}$  and therefore  $f_1(x)$ is computable.

Conversely, suppose that the function $f_1(x)$ of \eqref{pr3}
is computable. Let $W$ be a  word of length $|W| \le x$. Then $
W \in \langle \langle \R \rangle \rangle$ if and only if there
exists a van Kampen diagram $\D_W$  over \eqref{pr3} such that
$\ph(\p \D_W) \equiv W$ and $| \D_W(1)| \le f_1(x)$. Note that
if $\Pi $ is a face (= closure of a 2-cell) in $\D_W$ then $|\p
\Pi | \le 2 | \D_W(1) |$, where $|\p \Pi |$ is the perimeter of
$ \Pi$. Since $\R$ is decidable, we can write down all words $R
\in \R$ with $|R| \le 2f_1(x)$ and hence we can construct all
possible diagrams $\D$ over \eqref{pr3} with $| \D(1)| \le
f_1(x)$ to determine whether $W \in \langle \langle \R \rangle
\rangle$. This proves that the set $\langle \langle \R \rangle
\rangle$ is decidable and the word problem is solvable for
\eqref{pr3}. \qed

\medskip
{\em Proof of Theorem~1.3.}  Suppose that $G = \langle a_1, \dots, a_m
\rangle$ is a group generated by elements $a_1, \dots, a_m$ and the
problem  to decide whether a given word over the alphabet $\A^{\pm 1} =
\{  a_1^{\pm 1}, \dots, a_m^{\pm 1} \}$,  represents the identity element
of $G$ is in {\bf NP}. Consider the presentation $G = \langle \A \, \| \,
\R \rangle$, where $\R$ consists of all nonempty cyclically reduced words
$R $ over  $\A^{\pm 1}$ that represent the identity element $G$. Then, by
the definitions, the problem to determine whether $U \in \R$ is in {\bf
NP} and the corresponding function $f_1(x)$ is bounded by the polynomial
$x$.

Conversely, suppose that $G$ is defined by a presentation $G = \langle \A
\, \| \,  \R \rangle$ such that the corresponding function $f_1(x)$ is
bounded by a polynomial $p(x)$ and the problem to decide whether a word
$U$ belongs to $\R$ is in {\bf NP}. For definiteness, assume that this
problem can be solved in nondeterministic time bounded by a polynomial
$q(|U|)$. We need to show that the word problem for $G = \langle \A \, \|
\, \R \rangle$ is in {\bf NP}.

A word $U$ over  $\A^{\pm 1}$  represents the identity  element
$G$ if and only if there exists a van Kampen diagram $\D_U$
such that  $\ph(\p \D_U ) \equiv U$ and $|\D_U(1)| \le f_1(|U|)
\le p(|U|)$. We also remark that $|\D_U(2)| \le 2 |\D_U(1)| \le
2p(|U|)$ and, for every face $\Pi$ of $\D_U$, we have $|\p \Pi
| \le 2 |\D_U(1)| \le 2p(|U|)$. Therefore, a van Kampen diagram
$\D$, where $|\D(1)| \le p(x)$, $|\D(2)| \le 2p(x)$ and $|\p
\Pi | \le 2p(x)$ for every face $\Pi$ in $\D$,  can be used as
a certificate to verify that a given word $U$ with $|U| \le x$
represents the  identity  element of $G$. This can be done in
time bounded by the polynomial
\begin{gather*}
 |\D(2) | \cdot  q(|\p \Pi |) + |\p \D|
\le  2p(x)\cdot  q(2p(x)) +x \ ,
\end{gather*}
where   $q(|\p \Pi |)$ estimates the time needed to verify that the word
$\ph(\p \Pi)$ is in $\R$ and $|\p \D|$   bounds  the time needed to check
that $\ph(\p \D) \equiv U$.

We remark that the above argument is retained with ${\bf
PSPACE}$ in place of ${\bf NP}$ (recall that ${\bf NPSPACE} =
{\bf PSPACE}$, see \cite{Pap}). \qed

\medskip

{\em Proof of Theorem~1.4.}  Part (a). \ Let $\D$ be a van
Kampen diagram over \eqref{pr3}.  Assuming that $|\D(1)|>0$, it
is easy to see that $|\D(0)| \le 2 |\D(1)|$ and $|\D(2)| \le 2
|\D(1)|$. In the notation   \eqref{Lj}, these inequalities mean
that $L_0(W) \le 2 L_1(W)$ and $L_2(W) \le 2 L_1(W)$.
Maximizing over all $W \in \langle
 \langle \R \rangle  \rangle$ with $|W| \le x$,  we obtain
$ f_0(x)  \le 2f_1(x)$ and $f_2(x)  \le 2f_1(x)$,
 as required.
\medskip

Part (b). \  Let $\D$ be a van Kampen diagram over \eqref{pr3},
$\ph(\p \D) \equiv U$ and $|\D(1)|$ is minimal over all
diagrams $\D_1$ such that $\ph(\p \D_1) \equiv U$ and
$|\D_1(0)| = |\D(0)|$. To simplify the notation, let $V {}=
|\D(0)|$, $E {}= |\D(1)|$, $F {}= |\D(2)|$. Let  $F_2$ and
$F_3$ denote the numbers of faces in $\D$ that have 2 and $\ge
3$, respectively, edges in their boundaries. Recall that $\R$
has no words of length 1, whence $F= F_2 +F_3$. We also
consider a planar 2-complex $\D'$ obtained from $\D$ by
identifying $e$, $f$ for every pair of (oriented) edges $e$,
$f$ such that $ef^{-1}$ is the boundary cycle of a face of
$\D$. Let $V'$, $E'$, $F'$, $F'_2, F'_3$ be defined for $\D'$
in the same manner as the numbers $V$, $E$, $F$, $F_2, F_3$
were defined for $\D$. Clearly, $F'_2 = 0$ and
\begin{equation}\label{es1}
E = E' + F_2  \ .
\end{equation}

Suppose that $e_1$, $e_2$ are edges in $\D$ such that $(e_1)_-
= (e_2)_-$ and $(e_1)_+ = (e_2)_+$ (perhaps, $(e_1)_- =
(e_2)_+$), where $e_-$ denotes the initial vertex of an
oriented  edge $e$ and $e_+$  is the terminal vertex of $e$.
Also, assume that $\ph(e_1) = \ph(e_2)$ and the subdiagram,
bounded by the closed path $e_1 (e_2)^{-1}$, has no other
vertices  than $ (e_1)_-$, $(e_2)_+$. If $e_1 \ne e_2$, then
one could make a surgery on $\D$ that would identify $e_1$,
$e_2$, decrease $| \D(1)|$ and preserve both $| \D(0)|$ and
$\ph(\p \D)$. Since $\D$ is minimal relative to $|\D(1)|$, it
follows that the inequality $e_1 \ne e_2$ is impossible and so
$e_1 = e_2$ for such a pair $e_1$, $e_2$. This remark implies
that
\begin{equation}\label{es2}
F_2 \le 2m E'  \ ,
\end{equation}
where $m$ is the number of letters in $\A$. By the Euler
formula applied to $\D'$, we have $V' - E' +F' = 1$ or $V - E'
+F_3 = 1$, because $V' = V$, $F'_3 = F_3$ and $F'_2=0$. Note
that $F'_3 \le \tfrac{2E'}{3}$, hence $V - E' + \tfrac{2E'}{3}
\ge  1$ and $\tfrac{E'}{3} \le V$. By \eqref{es1}--\eqref{es2},
we obtain
\begin{equation}\label{D34}
|\D(1) | = E = E' +F_2 \le (1+2m)E' \le  3 (1+2m)V = 3
(1+2m)|\D(0) | \ .
\end{equation}
This estimate, together with the minimality of $\D$, implies that
$$
f_1(x) \le 3 (1+2m) f_0(x) \ .
$$
Hence, $f_1(x) \preccurlyeq f_0(x)$ which,
together with proven part (a), proves the equivalence of functions
$f_0(x)$ and $f_1(x)$. \qed

\medskip

{\em Proof of Theorem~1.5.} Recall that if $e$ is an oriented
edge of a diagram $\D$, then $e_-$ denotes the initial vertex
of $e$ and $e_+$  is the terminal vertex of $e$.  We will say
that a van Kampen diagram  $\D$ over \eqref{pr3} is {\em
$1$-regular} if $\D$ has the following property. If $e_1$,
$e_2$ are edges in $\D$ such that  $(e_1)_- = (e_2)_-$ and
$(e_1)_+ = (e_2)_+$ (perhaps, $(e_1)_- =  (e_2)_+$), $\ph(e_1)
= \ph(e_2)$ and the subdiagram $\D_e$ with $\p \D_e = e_1
e_2^{-1}$  has no  other vertices than  $(e_1)_-$, $(e_2)_+$,
then $e_1 = e_2$. Observe  that if $\D$ is not $1$-regular and
$e_1$, $e_2$ are the edges in $\D$ that violate the above
property, then one could make a surgery on $\D$, called {\em
$1$-reduction}, that would take $\D_e$ out of $\D$ and identify
$e_1$, $e_2$. Note that a $1$-reduction  preserves both $\ph(\p
\D)$ and $|\D(0)|$ and it  decreases $| \D(1)|$. Hence,
application of finitely many  1-reductions to $\D$ will yield a
1-regular  diagram $\D_1 $ such that  $\ph(\p \D_1)\equiv
\ph(\p \D)$ and  $| \D_1(0)|= | \D(0)|$. In particular, without
loss of generality,  we may assume that if $U \in \langle
 \langle \R \rangle  \rangle$ and  $\D$ is a  van Kampen
diagram over \eqref{pr3} such that $\ph(\p \D) \equiv U$ and
$ | \D(0)| = L_0(U)$, see \eqref{Lj}, then $\D$ is $1$-regular.

From now on assume that $|R| > 1$ for every $R \in \R$. Repeating the
arguments of the proof of Theorem~1.4(b) aimed to prove the inequality
\eqref{D34}, we can analogously show that if $\D$ is a $1$-regular van
Kampen diagram over \eqref{pr3}, then
\begin{equation}\label{D34a}
|\D(1) |  \le  3 (1+2m)|\D(0)| \ .
\end{equation}

Now suppose that the word problem is solvable for \eqref{pr3}
and $\R$ is decidable. Then $\langle
 \langle \R \rangle  \rangle$ is also a decidable subset of
 $F(\A)$ and, as in the proof of Theorem~1.2, for every $U \in
\langle \langle \R \rangle  \rangle$ we can effectively
construct a van Kampen diagram  $\D_U$  over  \eqref{pr3} with
$\ph(\p \D_U) \equiv U$. Applying 1-reductions to  $\D_U$   if
necessary, we may assume that $\D_U$  is 1-regular. By
inequality  \eqref{D34a}, $|\D_U(1) |  \le  3 (1+2m)|\D_U(0)
|$. Now let $\D_{U, 0}$ be a van Kampen diagram over
\eqref{pr3} such that $\ph(\p \D_{U, 0}) \equiv U$ and
 $| \D_{U,0}(0)| = L_0(U)$. As was pointed out above,
we may assume that $\D_{U, 0}$ is 1-regular, hence, by
\eqref{D34a},  $|\D_{U, 0}(1) |  \le  3 (1+2m)|\D_{U, 0}(0)|
\le   3 (1+2m)|\D_{U}(0)|$. This means that by checking all
diagrams $\D$ that satisfy $ |\D(1) |  \le  3
(1+2m)|\D_{U}(0)|$ we can compute the number $L_0(U)$. Thus
$f_0(x) = \max \{ L_0(U) \mid U \in  \langle  \langle \R
\rangle \rangle , \ |U| \le x \}$ is also computable.

Conversely, suppose that the function $f_0(x)$ is computable
for \eqref{pr3} and let  $W$ be a  word with $|W | \le x$. Then
$W \in \langle \langle \R \rangle  \rangle$ if and only if
there exists a 1-regular van Kampen diagram  $\D_W$  over
\eqref{pr3}  such that  $\ph(\p \D_W) \equiv W$ and $|\D_{W}(0)
| \le f_0(x)$. In view of inequality  \eqref{D34a},
 \begin{equation}\label{D34b}
|\D_W(1) |  \le  3 (1+2m)|\D_W(0)| \le  3 (1+2m)f_0(x) \ .
\end{equation}
Hence, as in the proof of Theorem~1.2, we can construct all
possible diagrams $\D$ over \eqref{pr3}  with $|\D(1) |  \le  3
(1+2m)f_0(x)$ and determine whether or not  $W \in \langle
\langle \R \rangle  \rangle$. This proves that the word problem
is solvable for \eqref{pr3}.

Finally, by Theorem~1.2, the Dehn 1-function $f_1(x)$ of
\eqref{pr3} is computable if and only if the word problem for
\eqref{pr3} is solvable and, as was shown above (when $\forall
R \in \R$ $|R| > 1$),  the Dehn 0-function $f_0(x)$ of
\eqref{pr3} is computable if and only if the word problem  for
\eqref{pr3} is solvable. This shows that $f_0(x)$ is computable
if and only if so  is $f_1(x)$.
 \qed

\medskip

{\em Proof of Theorem~1.6.}  Part (a). \ Let $\R' = (\R
\setminus \s) \cup \U$ and $\s = \{ S_1, \dots, S_k\}$,
 $\U = \{ U_1, \dots, U_\ell \}$. Since $\U \subset \langle
 \langle \R \rangle  \rangle$, there are van Kampen diagrams
 $\D_i$ over the presentation $  \langle
 \A \, \| \,  \R'  \rangle$ such that
 $\ph(\p \D_i) \equiv S_i$, $i = 1, \dots, k$. Denote
 \begin{gather*}
 M_j {}=  \max \{ | \D_i(j) | \mid  i = 1, \dots, k \} \ ,
 \quad j = 0,1,2 \  .
 \end{gather*}
Consider a van Kampen diagram  $\D$ over \eqref{pr3} and let $\Pi$ be a
face in $\D$ with $\ph(\p \Pi) \equiv S_i^\e$, $\e = \pm 1$. We replace
$\Pi$ in $\D$ by a copy of $\D_i$ if $\e =1$ or by a  mirror copy of
$\D_i$ if $\e =-1$. Doing this for all faces $\Pi$ in $\D$ with $\ph(\p
\Pi) \equiv S_i^\e$, where $i= 1, \dots, k$ and $\e = \pm 1$, results in
a diagram $\D'$ over the presentation $\langle \A \, \| \, \R' \rangle$.
Observe that
\begin{gather*}
 \max (|\D(0) |, |\D(2) |)  \le 2 |\D(1)| \ , \quad |\D'(0) |
\le  |\D(0) |+ M_0|\D(2) |  \ , \\    |\D'(1) | \le |\D(1) |+
M_1|\D(2) |  \ , \quad     |\D'(2) | \le |\D(2) | + M_2|\D(2) |
\  .
\end{gather*}
Hence,
\begin{gather*}
|\D'(0) | \le  2(1+ M_0) |\D(1)|   ,  \  |\D'(1) | \le (1+2
M_1)|\D(1) |   , \     |\D'(2) | \le  (1+ M_2)|\D(2) | \ .
\end{gather*}
It follows from these inequalities and the definitions that
\begin{equation*}
f'_0(x)  \le  2(1+ M_0)f_1(x) \ ,  \ f_1'(x) \le  (1+
2M_1)f_1(x) \ ,  \ f'_2(x) \le (1+ M_2)f_2(x) \ .
\end{equation*}
Therefore, $f'_0(x)  \preccurlyeq f_1(x)$, $f'_1(x)
\preccurlyeq f_1(x)$, $f_2'(x) \preccurlyeq f_2(x)$ which, in
view of  symmetry between $\R$ and $\R'$, imply the desired
relations.
\medskip

Part (b). \ Let a presentation $\langle \A' \, \| \, \R'
\rangle$ be obtained from \eqref{pr3} by a stabilization, $\A'
= \A \cup \{ b \}$, and $f'_j(x)$, $f_j(x)$, $j=0,1,2$, be
their corresponding Dehn functions. Note that if $\D$ is a
diagram over $\langle \A' \, \| \, \R' \rangle$ then every edge
$e$ of $\D$ with $\ph(e) = b^{\pm 1}$ lies on the boundary $\p
\D$ of $\D$. This remark and the definitions enable us to
conclude that $f_j(x) \le f'_j(x)$ and $f_j'(x) \le f_j(x) +x$,
$j=0,1,2$. These inequalities imply the required equivalence
$f_j'(x) \simeq f_j(x)$, $j=0,1,2$. \qed
\medskip

{\em Proof of Theorem~1.7.}  \ Let $\R$ be finite. Then  $M =
\max\{ |R| \mid R \in \R \}$ is also finite and $|\D(1)| \le
M|\D(2)| + |\p \D |/2$ for every diagram $\D$  over
\eqref{pr3}, where $|\p \D |$ is the perimeter of $\D$. Using
the notation of the definition \eqref{DFj}, we further have
$L_1(W) \le L_2(W) +|W|/2$. Maximizing over all $W$, where $W
\in \langle \langle \R \rangle \rangle$ and $|W| \le x$, we get
$f_1(x) \le M f_2(x) + x/2$ and so  $f_1(x) \preccurlyeq
f_2(x)$. This, together with Theorem~1.4(a), proves the
equivalence $f_1(x) \simeq f_2(x)$.

To prove the equivalence  $f_0(x) \simeq f_1(x)$, we will argue
by induction on the total length $\| \R \| = \sum_{R \in
\R}|R|$ of words in $\R$. The base step for $\| \R \| \le 1$ is
obvious. If all words in $\R$ have length $> 1$, then the
desired equivalence follows from Theorem~1.4(b). Without loss
of generality, we may assume that $\R$ contains a letter $b \in
\A$ and the set $\R^b {}= \R \setminus \{ b^{} \}$ is nonempty.

If $b, b^{-1}$ do not occur in words of $\R^b$, then we can
apply a stabilization to \eqref{pr3} and obtain  the
presentation $\langle \A \setminus \{ b \} \, \| \, \R^b
\rangle$.  By the induction hypothesis, its Dehn functions
$f_0^b(x)$, $f_1^b(x)$ are equivalent and, by Theorem~1.6(b),
$f_j^b(x)  \simeq f_j(x)$, $j=0,1$. Hence, $f_0(x)$, $f_1(x)$
are also equivalent, as required.

Now assume that $b$ (or $b^{-1}$) occurs in a relator $R_b \in
\R^b$. Let
\begin{equation}\label{pds}
R_b \equiv  U b U^{-1} S \ ,
\end{equation}
where $U, S \in F(\A)$ and $S$ is cyclically reduced. Clearly, $S \in
\langle\langle \R^b \rangle \rangle $. Set
$$
\R' {}= (\R \setminus \{ R_b \} ) \cup \{ S \}
$$
 and consider the presentation
\begin{equation}\label{npr}
\langle \, \A \ \| \ \R'  \, \rangle \ .
\end{equation}
Note  $\|\R' \| < \| \R \|$. Pick a cyclically reduced word $W
\in \langle\langle \R \rangle \rangle $ and let $\D$ be a van
Kampen diagram  over  \eqref{pr3} such that $\ph(\p \D) \equiv
W$ and $\D$ is minimal relative to $| \D(0) |$. If $\Pi$ is a
face in $\D$ and $\ph(\p \Pi) \equiv R_b^\e$, $\e =\pm 1$, then
we consider a subpath $p$ of $\p \Pi = p q$ whose label
$\ph(p)$ is the  subword $U b^\e U^{-1}$ of $R_b^\e$
distinguished in \eqref{pds}. Let $p = uev$, where $\ph(u)
\equiv \ph(v)^{-1} \equiv U$, $\ph(e) \equiv b^{\e}$. If the
initial vertex $e_-$ of $e$ were  different from its terminal
vertex $e_+$, then we could put two faces $\pi, \pi'$, with
$\ph(\p \pi) =  \ph(\p \pi')^{-1} = b^\e$, into $\D$, see
Figure~1, thus making $e_-$, $e_+$ merge and decreasing $|
\D(0) |$ by one. This contradiction to the minimality of $|
\D(0) |$ shows that $e_- = e_+$.

\unitlength .6mm \linethickness{0.4pt}
\ifx\plotpoint\undefined\newsavebox{\plotpoint}\fi 
\begin{picture}(96.00,38.00)(-36.00,-2.00)
\put(25.00,20.00){\line(1,0){15.00}}
\put(55.00,20.00){\line(1,0){13.00}}
\put(90.00,26.00){\circle{12}} \put(90.00,14.00){\circle{12}}
\put(25.00,20.00){\circle*{1.50}}
\put(40.00,20.00){\circle*{1.50}}
\put(90.00,20.00){\circle*{1.50}}
\put(32.00,20.00){\vector(1,0){2.75}}
\put(59.50,20.00){\vector(1,0){5.25}} \put(31.50,23.50){$e$}
\put(88.75,24.50){$\pi$} \put(88.50,12.75){$\pi'$}
\put(50.75,4.75){Figure~1}
\end{picture}

Let $u = u_1 u_2$,  $v = v_1 v_2$ be some factorizations of $u,
v$ such that $|u_2| = |v_1|$, where $|u_2|$ denotes  the length
of $u_2$, and $(u_2)_- = (v_1)_+$, see Figure~2. We pick such
factorizations so that $|u_2|$ is maximal (perhaps, $|u_2|=0$).
Since $(u_2)_- = (v_1)_+$ and $\ph(u_1) \equiv \ph(v_2)^{-1}$,
we can do the following surgery over $\D$. Take the subdiagram
bounded by the closed path $u_2 e v_1$ out of $\D$ and identify
the paths $u_1$ and $v_2^{-1}$, see Figure~2.

\unitlength 0.64mm \linethickness{0.4pt}
\ifx\plotpoint\undefined\newsavebox{\plotpoint}\fi 
\begin{picture}(162.75,54)(0.00,4.00)
\put(25.00,35.00){\line(1,0){55.00}}
\put(96.25,35.00){\line(1,0){14.75}}
\put(129.75,35.00){\line(1,0){33.00}}
\put(145.50,35.00){\line(0,1){15.00}}
\qbezier(53.00,35.00)(47.00,27.50)(53.00,22.00)
\qbezier(53.00,35.00)(59.00,27.50)(53.00,22.00)
\put(53.00,18.00){\circle{8.00}}
\put(53.00,35.00){\circle*{1.50}}
\put(53.00,22.00){\circle*{1.50}}
\put(70.00,35.00){\circle*{1.50}}
\put(36.00,35.00){\circle*{1.50}}
\put(145.50,35.00){\circle*{1.50}}
\put(145.50,50.00){\circle*{1.50}}
\put(64.00,35.00){\vector(-1,0){3.00}}
\put(46.25,35.00){\vector(-1,0){2.75}}
\put(28.25,35.00){\vector(-1,0){1.75}}
\put(136.25,35.00){\vector(-1,0){2.75}}
\put(145.50,41.50){\vector(0,1){2.75}}
\put(102.50,35.00){\vector(1,0){4.25}}
\put(56.00,28.50){\vector(0,-1){1.75}}
\put(50.00,27.25){\vector(0,1){2.00}}
\put(53.50,14.00){\vector(-1,0){1.50}}
 \put(92.00,12){Figure 2}
 \put(61.00,31.00){$u_1$}  \put(57.500,25.90){$u_2$}
\put(53.00,9.50){$e$} \put(44.5,26.){$v_1$}
\put(43.,31.0){$v_2$} \put(27.00,30.){$q$}
\put(68.50,20.75){$\Pi$} \put(144.50,25.50){$\Pi'$}
\put(148.75,42.75){$u_1 = v_2^{-1}$} \put(133.75,30.){$q$}
\put(123.75,40.){$\p \Pi' = q$}  \put(30,44){$\p \Pi = u_1 u_2
e v_1 v_2 q $}

\end{picture}

\noindent By doing this, we turn the face $\Pi$ with $\ph(\p
\Pi) = R_b^\e$ into a face $\Pi'$ with $\ph(\p \Pi') = S^\e$
and do not increase $|\D(0) |$. Iterating such surgeries for
all faces $\Pi$ as above, we will obtain a diagram $\D'$ over
the presentation \eqref{npr} such that $\ph(\p \D') \equiv W$
and
\begin{equation}\label{dp0}
|\D'(0) | \le |\D(0) | \ .
\end{equation}
It follows from the choice of $\D$  that $L_0(W) = |\D(0) |$.
Hence, referring to the inequality \eqref{dp0} and the
definition \eqref{DFj}, we obtain  that
\begin{equation}\label{lst1}
f_0'(x) \le    f_0(x)  \ ,
\end{equation}
where $f_0'(x)$ is the Dehn 0-function of  the presentation
\eqref{npr}. Since the presentation \eqref{npr} is obtained
from   $\langle \, \A \, \| \, \R  \, \rangle$ by a
$T$-transformation which replaces $R_b$ by $S$, it follows from
Theorem~1.6(a) that $f_1'(x) \simeq f_1(x)$. In view of $\| \R'
\| < \| \R \|$, the induction hypothesis applies to \eqref{npr}
and yields that $f'_0(x) \simeq f'_1(x)$. Hence,
\begin{equation}\label{lst2}
f'_0(x) \simeq f'_1(x) \simeq f_1(x)   \ .
\end{equation}
By Theorem~1.4(a), $f_0(x) \preccurlyeq f_1(x)$ which, in view
of \eqref{lst2}, means that $f_0(x) \preccurlyeq f_0'(x)$.
This, together with \eqref{lst1}, shows that $f_0'(x) \simeq
f_0(x)$ and, by \eqref{lst2}, we finally have $f_0(x) \simeq
f_1(x)$. Theorem~1.7 is proved. \qed

\medskip

Let us give  examples of group presentations  for which the
pairs $( f_1(x), f_2(x))$,  $( f_0(x), f_1(x))$, $( f_0(x),
f_2(x))$ contain nonequivalent functions.
\medskip

\noindent {\bf Example 2.1}. \ The presentation
 \begin{equation*}
 \langle \ a, b  \  \| \
a^iba^{-i}b^{-1}, \  i \in \mathbb N \,  \rangle
 \end{equation*}
defines a free abelian group of rank 2.  It is easy to check
that  $f_1(x)\simeq x^2$ and $f_2(x) \simeq x$. Hence, $f_1(x)
\not\simeq f_2(x)$ for this presentation.
\medskip

\noindent
{\bf Example 2.2}. Consider the presentation
\begin{equation*}
 \langle \ a, b, c  \  \| \ c, \
R_i c^{\ell_i}, \  i \in \mathbb N \,   \rangle \ ,
 \end{equation*}
where $R_i$ are words over positive alphabet $\{ a, b\}$ that
satisfy the small cancelation condition $C'(\lambda)$, $0 <
\lambda < 1/6$, see \cite{LS},  $| R_i | \to \infty$ as $i \to
\infty$, and  $\ell_i = 2^{| R_i| }$. It is not difficult to
verify that $f_0(x) \simeq x$ and $\min (f_1(x), f_2(x)) \ge
2^x$ for every $x$ for which there exists $R_i$ with $| R_i| =
x$. Thus, the functions $f_0(x)$, $f_1(x)$ are not equivalent
and $f_0(x), f_2(x)$ are not equivalent either.
\medskip

It is of interest to point out that, when proving Theorem~1.2
(resp. Theorem~1.5), we actually show  that the function
$L_1(W)$ (resp. $L_0(W)$), where $W \in \langle \langle \R
\rangle \rangle$, see \eqref{Lj}, is computable if  the word
problem is solvable for \eqref{pr3}. The following example, due
to Jockush and Kapovich, gives an indication that
Problem~1.1(a) for $j =2$  might have a negative solution.
\medskip

\noindent
 {\bf Example 2.3}. Consider the presentation
\begin{equation*}
 \langle \ a, b \  \| \ a^i, \  a^i b^{k_i},  \ \
 i \in \mathbb N \,   \rangle \ ,
 \end{equation*}
where $\mathbb K = \{ k_1, k_2, \dots \}$ is a recursively enumerable but
not recursive subset of $\mathbb N$ with the indicated enumeration and
$k_1 =1$. It is clear that the relator set is decidable and this
presentation defines the trivial group, hence the word problem is
solvable. On the other hand, it is easy to verify that $L_2(b^k) = 2$,
where $k \in \mathbb N$,  if and only if $k \in \mathbb K$. Since
$\mathbb K$ is not recursive , it follows that the function $L_2(W)$,
where $W$ is a word over $\{ a^{\pm 1}, b^{\pm 1} \}$, is not computable.
It remains to be seen whether this idea would lead to  a counterexample
to Problem~1.1(a) for $j=2$.
\medskip

The following example that gives a negative solution to
Problem~1.1(b) for $j=2$ is due to an anonymous referee.

\medskip

\noindent {\bf Example 2.4}. Let $\langle  \A_0 \,  \| \, \R_0 \rangle$
be a finite presentation with unsolvable word problem and assume that
$\R_0$ contains a letter $a$ of $\A_0$. Consider a new letter $t$, $t
\not\in \A_0^{\pm 1}$, denote $ \A =  \A_0 \cup \{ t \}$ and let $\N_t$
be  the set of all nonempty cyclically reduced words over $\A^{\pm 1}$
that are in the normal closure of $t$ in $F(\A)$. Observe that $ W \in
\N_t$ if and only if $W$ is nonempty, cyclically reduced and $\pi_t(W)
=1$, where $\pi_t : F(\A) \to F(\A_0)$ is the projection homomorphism
that erases all occurrences of $t^{\pm 1}$. For every $2k$-tuple $(R_1,
X_1, \dots, R_k, X_k)$, where $k \ge 1$, $R_i \in \R_0^{\pm 1}$ and $X_i$
are reduced words over $\A_0^{\pm 1}$ for  $i =1, \dots, k$, we consider
the word
$$
V(R_1, X_1, \dots, R_k, X_k)=  X_1 t R_1 t X_1^{-1} t  X_2 t
R_2 t X_2^{-1} t \dots X_k t R_k t X_k^{-1} t \ .
$$
 Let the set $\V_t$ contain
the words $V(R_1, X_1, \dots, R_k, X_k)$ for all possible
$2k$-tuples $(R_1, X_1, \dots, R_k, X_k)$, as described above.
Now we define the presentation
\begin{equation}\label{e24}
 \langle \ \A  \  \| \  \R =  \N_t\cup \V_t  \,   \rangle \ .
 \end{equation}
It is clear that $\R$ is a decidable set.
 Observe that
 $$
 \pi_t(V(R_1, X_1, \dots, R_k,
X_k) ) = X_1 R_1 X_1^{-1}  \dots   X_k R_k X_k^{-1} \ , \quad
 \pi_t( \N_t ) = \{ 1 \} \ ,  \quad t \in \N_t \ .
 $$
 Therefore, a word
$W$ over $\A^{\pm 1}$ is in $\langle \langle   \R  \rangle
\rangle$ if and only if $\pi_t(W)$ is in the normal closure
$\langle \langle   \R_0 \rangle \rangle^{F(\A_0)}$    of the
set $\R_0$ in $F(\A_0)$. In particular, the word problem is
unsolvable for the presentation \eqref{e24}.

Now assume that $W \in \langle \langle   \R \rangle \rangle$.
Then $\pi_t(W) \in  \langle \langle   \R_0 \rangle
\rangle^{F(\A_0)}$ and hence there exists a suitable word $V
=V(R_1, X_1, \dots, R_k, X_k)$ in $\R$ such that $\pi_t(V) =
\pi_t(W)$. Since $\pi_t(W V^{-1}) = 1$, it follows that either
$W V^{-1} =1$ in $F(\A)$ or $W V^{-1}$ is conjugate in $F(\A)$
to a word in $\N_t$. Writing $W$ in the form  $W = (W V^{-1})
V$, we see that $L_2(W) \le 2$, where $L_2(W)$ is defined by
means of \eqref{e24}. Furthermore, it follows from the
definitions that $a   \not\in \R^{\pm 1}$ which implies that
$L_2(a) = 2$. Thus, the Dehn 2-function $f_2(x)$ of the
presentation \eqref{e24} is identically equal to 2 and hence is
computable, whereas the word problem for \eqref{e24} is
unsolvable.

\medskip

We remark that it is also interesting to state Problem~1.1  requiring, in
addition, that  $\R$ be minimal. The presentations of Examples~2.3--2.4
are far from being minimal and potential counterexamples (if they exist)
to the analog of Problem~1.1  with minimal $\R$  could be more difficult
to construct.

\medskip

As was suggested by the referee, it is of interest to consider an ``upper
bound" form of Problem~1.1 in which the computability of the Dehn
function $f_j(x)$ is replaced by  the computability of an upper bound of
$f_j(x)$. Observe that the upper bound version of Problem~1.1(a) has a
straightforward positive solution  for both $j=0, 2$. However, even this
relaxed  version of Problem~1.1 still has a negative solution for $j =2$,
as follows from Example~2.4.

The referee also  pointed out that our idea of the Dehn
1-function $f_1$ is analogous to the concept of ``derivation
work" introduced by Birget \cite{Bir} in 1998 for semigroup and
group presentations and, in fact, Birget's ``derivation work"
function is equivalent to $f_1(x)$ in terms of the linear
equivalence $\simeq$. Furthermore,  the left-to-right direction
of Theorem~1.3 is a consequence of Proposition~3.3 and
Corollary~3.4 of \cite{Bir} where a much stronger constraint
for $\R$ is obtained: the membership problem for $\R$ is not
only in {\bf NP} but $\R$ can be chosen to be the intersection
of two deterministic context-free languages which, in
particular, implies that the membership problem for $\R$ is
solvable in deterministic linear time.

\section{Proofs of Theorem~1.8 and  Corollary~1.9}

As in Introduction, let $\Ga $ denote the group introduced in \cite{G80}
and defined by presentation \eqref{pr1}. This group $\Ga$ turned out to
possess many interesting  properties: $\Ga$ is an infinite 3-generator
2-group all of whose nontrivial quotients are finite \cite{G80, G00},
$\Ga$ has bounded width with respect to the lower central series, $\Ga$
is of intermediate (between polynomial and exponential) growth
\cite{G84},  $\Ga$ is amenable but not elementary amenable \cite{G84}
etc. A detailed discussion of  properties of the  group $\Ga$ can be
found in  \cite{G84, G00, G05, harpe}. This  group $\Gamma$ can be
defined in several different ways but, in this article, we will only use
the  definition of  $\Gamma$    by means of the presentation
\eqref{pr1}.

\medskip

Let $\Ga(0)$ denote the free group $F(a,b,c,d)$ in the alphabet
$\{ a,b,c,d\}$ and $\R(0)$ be the empty set. Consider the group
presentation
\begin{equation}\label{ga1}
\Ga(1)   = \langle \ a, b, c, d \ \| \ a^2, \ b^2,\ c^2,\ d^2,
\ bcd, \ (ad)^4  \ \rangle \  = \langle \ a,b,c,d{{}} \ \| \
\R(1) \ \rangle  \ .
\end{equation}
For $i \ge 2$, we define
\begin{align} \notag
\Ga(i)   = \langle & \ a,b,c,d \ \| \ a^2, \ b^2,\ c^2,\ d^2, \
bcd, \ (ad)^4  ,  \  \sigma^j((ad)^4)  ,  \\ &
\sigma^{j-1}((adacac)^4) ,   \ j=1, \dots, i-1 \rangle \ = \
\langle \  a,b,c,d{{}} \ \| \ \R(i) \ \rangle \ . \label{gaii}
\end{align}

For every $i \in \mathbb N^* = \mathbb N \cup \{ 0 \} \cup \{
\infty \}$, let $\N(i)$ denote the normal closure
$\langle\langle \R(i) \rangle\rangle$ of the relator set
$\R(i)$ in $F(a,b,c,d)$. It follows from  the Lysenok theorem
\cite{Lys85} that $\Ga = F(a,b,c,d)/\N(\infty)$. Note that,
whenever it is not ambiguous, we do not distinguish between a
group and its presentation. Unifying the foregoing notation
\eqref{pr1}, \eqref{ga1}, \eqref{gaii},  we can define the
presentation
\begin{equation*}
\Ga(i) = \langle \  a,b,c,d{{}} \ \| \ \R(i) \ \rangle
\end{equation*}
for every $i \in \mathbb N^*$.

Consider the subgroup $\HH(i)$ of   $\Ga(i)$  that is generated
by the  images of the words
\begin{equation}\label{abc}
b,\ c, \ d, \  aba, \ aca, \ ada, \ a^2 \ .
\end{equation}
It is easy to see that $\HH(i)$ is of index 2 in $\Ga(i)$.

Observe that the words \eqref{abc} are free generators of the
free group $\HH(0)$. Define a homomorphism
$$
\psi_0 : \HH(0) \mapsto \Ga(0) \times \Ga(0)
$$
by setting
\begin{align}\label{psi0}
\psi_0 {}= \begin{cases}
  b    \mapsto (a,c) , \\ c  \mapsto (a,d), \\ d  \mapsto
 (1,b) , \\   aba  \mapsto  (c,a) , \\ aca
 \mapsto  (d,a), \\ ada   \mapsto  (b,1) \\ a^2  \mapsto
(1,1) \ .
 \end{cases}
\end{align}

Referring to the definition of $\sigma$, see \eqref{sig}, we
observe that $\sigma(\Ga(0)) \subseteq \HH(0)$. Hence, we can
compose
$$\psi_0 \cdot \sigma : \Ga(0) \mapsto
\Ga(0) \times \Ga(0) \ .
$$
Computing $\psi_0 \cdot \sigma$, we obtain
\begin{align}\label{psisig}
\psi_0 \cdot \sigma = \begin{cases} a    \mapsto (d,a) , \\  b
\mapsto (1,b) , \\ c  \mapsto (a,c), \\ d  \mapsto
 (a,d) \ .
 \end{cases}
\end{align}

It will be convenient to partition the relator set $\R(\infty)$
of the presentation $\Ga(\infty)$, or  \eqref{pr1}, as follows.
Let
$$
\s(0) =  \{ a^2, \ b^2, \ c^2, \ d^2, \ bcd , \ (ad)^4 , \
(adacac)^4 \}
$$
and, for $i \ge 1$, we set $\s(i) = \{ \sigma^i((ad)^4),
\sigma^i((adacac)^4) \}$. Clearly,
$$
\R(\infty) =\bigcup_{i=0}^\infty \s(i)
$$
is a partition of $\R(\infty)$.  The {\em height} $h(R)$ of a
relator $R \in \R(\infty)$ is defined to be  $j\ge 0$ if $R \in
\s(j)$.

\begin{Lem31}\label{Lem1}
For every $i \in \mathbb N^*$, the map $\psi_0$, defined by
\eqref{psi0}, extends  to a homomorphism
\begin{align}\label{psii}
 \psi_i : \HH(i) \mapsto \Ga(i) \times \Ga(i)
\end{align}
such that   $\ker \psi_\infty = \{ 1\}$ and $\ker \psi_1 =
\N(2) / \N(1)$. In particular, $\ker \psi_1$ is a normal
subgroup of $\Ga(1)$ generated by the images of $\sigma((ad)^4)
= (ac)^8$ and $(adacac)^4$.
\end{Lem31}

\begin{proof}

Referring to the classic Reidemeister--Schreier rewriting
process, see \cite{LS} or \cite{MKS}, we conclude that, to
prove that $\psi_0$ extends  to a homomorphism \eqref{psii}, it
suffices to verify that
\begin{equation}\label{inc}
\psi_0(\R(i)) \cup  \psi_0(a\R(i)a) \subset \N(i) \times \N(i)
\ .
\end{equation}
It is easy to check that
$$
\psi_0(\s(0) \cup a \s(0) a ) \subset \N(1) \times \N(1) \ .
$$
In fact, to show this inclusion we do not even need to use the
inclusion $(ad)^4 \in \N(1)$. Since $\R(1) \subset \s(0)$, it
follows that $\psi_0(\R(1) \cup  a\R(1)a) \subset \N(1) \times
\N(1)$.

Suppose $U \in \R(i)$, where $i \ge 2$, and $ U \not\in \s(0)$.
It follows from the definitions that  $U = \sigma(V^4)$, where
$V^4$ is the $\sigma^j$-image of $(ad)^4$ or $(adacac)^4$,
where $0 \le j \le i-2$, and  so $V^4 \in \R(i-2)$. Hence,
according to \eqref{psisig},
\begin{equation}\label{psiU}
 \psi_0(U) =  \psi_0(\sigma(V))^4 = (T^4, V^4) \ ,
\end{equation}
where $T$ is a word over the alphabet $\{ a, d\}$. Since $a^2,
d^2, (ad)^4 \in \N(1)$ and $V^4 \in \R(i-2)$, it follows that
$T^4 \in \N(1)$ and $\psi_0(U) \in \N(i) \times \N(i)$, as
desired.

It follows from the definition of $\psi_0$, see \eqref{psi0},
that, switching $g \leftrightarrow a g a$, where $g \in \{
b,c,d\}$, results in switching the first and second components
of $\psi_0(e)$. This remark, together with \eqref{psiU} and the
equality $\psi_0(a^2) = (1,1)$,  implies that
\begin{equation}\label{psiaUa}
 \psi_0(aUa) =
\psi_0(a\sigma(V^4)a) = (V^4 , T^4  ) \ .
\end{equation}
 Hence, we also have $\psi_0(aUa) \in \N(i) \times
\N(i)$ and the inclusion \eqref{inc} is proved.

 The equality $\ker \psi_\infty = \{ 1\}$ was
observed in the original article \cite{G80}, see also
\cite{harpe}, and the equality $\ker \psi_1 = \N(2) / \N(1)$ is
shown in \cite{GH}.
\end{proof}

\begin{Lem32}\label{Lem2} Let $U \in \N(\infty)$ and $\psi_0(U) = (U_0,
U_1)$. Then, in the free group $F(a,b,c,d)$,  $U = a
\sigma(U_0) a \sigma(U_1)V$, where $V \in  \N(2)$.
\end{Lem32}

\begin{proof} By Lemma~3.1, we have $U_0, U_1 \in \N(\infty)$
and, in view of $\sigma(\N(\infty))  \subset  \N(\infty)$, it
follows from \eqref{psisig} that
$$
 \psi_0(\sigma(U_1)) = (T_1,
U_1) \ ,
 $$
 where $T_1$ is a word over $\{ a, d\}$, and hence
$T_1 \in \N(1)$. As in the proof of Lemma~3.1, in view of the
symmetry of the definition of $\psi_0$, see \eqref{psiaUa},
relative to the switch $g \leftrightarrow a g a$, $g \in \{
b,c,d\}$, we analogously   obtain
$$
 \psi_0(a\sigma(U_0)a) = (U_0, T_2) \ ,
 $$
where $T_2  \in \N(1)$. Therefore,
$$
 \psi_0(a\sigma(U_0)a \sigma(U_1) ) = (U_0  T_1, T_2 U_1) \
 .
 $$
Since $\psi_0 (U) = (U_0, U_1)$, it follows from Lemma~3.1 that
$$
 U^{-1} a\sigma(U_0)a \sigma(U_1) \in \ker \psi_1 \cdot \N(1)  =
 \N(2) \ ,
$$
as required.
\end{proof}

Observe that if $R \in \s(i)$, $i \ge 1$, then $\sigma(R) \in
\s(i+1)$. On the other hand, if $R \in \s(0)$, then either
$\sigma(R) \in \s(1)$ or, up to a cyclic permutation (when $R =
bcd$), $\sigma(R) \in \s(0)$ except for the case $R = a^2$.
Since  $\sigma(a^2) \not\in \R(\infty)$, we wish to extend the
sets $\s(i)$, $i \ge 1$, by adding $\sigma^i(a^2)$. Set
$\s^*(0) = \s(0)$ and define
$$
\s^*(i) = \s(i) \cup \{ \sigma^i(a^2)\} \quad \ \text{for }
\quad \  i \ge 1 \ .
$$

Note that the map $\sigma$, see \eqref{sig}, extends to a
homomorphism $\Ga(\frac 12 ) \mapsto \Ga(\frac 12 )$ of the
free product $\Ga(\frac 12 ) = \langle \ a, b, c, d \ \| \ a^2,
b^2, c^2, d^2 \ \rangle$ of four groups of order 2. Hence,
$\sigma^j(a^2) \in \N(1)$ for every $j \ge 0$. In particular,
the addition of the relator $\sigma^j(a^2)$ to $\s(j)$,   $j
\ge 1$, does not change the normal closure of $\s(j)$. Hence,
letting $\R^*(\infty) = \bigcup_{j=0}^\infty \s^*(j)$, we
obtain another presentation for the group $\Ga$
\begin{equation}\label{ga*}
\Ga^*(\infty) = \langle \  a,b,c,d{{}} \ \| \ \R^*(\infty) =
\bigcup_{j=0}^\infty \s^*(j)  \ \rangle \ ,
\end{equation}
because $\langle\langle \R^*(\infty) \rangle\rangle= \N(\infty)$.

As above, if $R \in \s^*(j)$, then we say that $R$ has the {\em
height} $h(R) =j$.

Assume that $W $ is a word in $\N(\infty)$. As in \eqref{FD},
consider a product for $W$ of the form
\begin{equation}\label{DF2}
W =  \prod_{j=1}^L  X_j  R_j^{\e_j}  X_j^{-1}   \ ,
\end{equation}
where $X_j \in F(a,b,c,d)$, $R_j \in \R^*(\infty)$, and $\e_j =
\pm 1$. An {\em $h^*$-tuple} $\tau^W = ( \tau_0^W, \tau_1^W,
\dots)$ of a word $W \in \N(\infty)$ is defined by means of a
product  \eqref{DF2} for $W$  so that $\tau_i^W$ is the number
of factors $X_j  R_j^{\e_j}  T_j^{-1}$ in \eqref{DF2} with
$h(R_j) = i$. Clearly,
\begin{equation*}
f^*_{2,{\Gamma}}(x) \le \max_{|W| \le x}  \sum_{j=0}^\infty \tau_j^W  \ ,
\end{equation*}
where $f^*_{j, {\Gamma}}(x)$ denotes  the Dehn $j$-function of
the presentation $\Ga^*(\infty)$, see \eqref{ga*}, $j =0,1,2$.

Putting the more  restrictive inclusion $R_j \in \R(\infty)$ in
the foregoing definition, we analogously define an {\em
$h$-tuple} $\bar \tau^W$ of a word $W \in \N(\infty)$ with
respect to the smaller relator set $\R(\infty) \subset
\R^*(\infty)$.

\begin{Lem33}\label{Lem3} Let $W \in \N(\infty)$ and $|W| \le x$.
Then $W$ possesses an  $h^*$-tuple $\tau^W = ( \tau_0^W,
\tau_1^W, \dots)$ such that $\tau_i^W = 0$ if $i > \log_2 x$,
$\tau_0^W \le Cx^2$, and $ \tau_i^W \le \frac {C}{2^{i-1}} x^2$
where $1 \le i \le \log_2 x$ and $C >1 $ is a constant
(independent of $W$). In particular, the Dehn  functions
$f^*_{1,{\Gamma}}(x)$, $f^*_{2,{\Gamma}}(x)$ of the
presentation $\langle \ a,b,c,d{{}} \ \| \ \R^*(\infty) \
\rangle$  of $\Ga$ satisfy the  following inequalities
\begin{align}\label{2inq1}
f^*_{2, \Ga}(x)  & \le \max_{|W| \le x}  \sum_{j=0}^\infty \tau_j^W \le
3Cx^2 \ , \\ \label{2inq2}  f^*_{1, \Ga}(x) & \le
 \max_{|W| \le x} \sum_{j=0}^\infty 3\cdot 2^{j+3}\tau_j^W \le 50  Cx^2\log_2 x \
.
\end{align}
\end{Lem33}

\begin{proof} By induction on $x \ge 1$, where $W$ is a
 word in $\N(\infty)$ with $|W| \le x$, we will be proving that
$W$ has a desired $h^*$-tuple $ \tau^W = ( \tau_0^W, \tau_1^W,
\dots)$. The base of induction for $x = 1$ is trivial and we
assume that $x \ge 2$.

First we note that, using at most $3x$ copies of relators
$a^{2}$, $b^{2}$, $c^{2}$, $d^{2}$, $bcd^{}$ (when ``using'' we
allow cyclic permutations and inversions), we can turn the word
$W$ into a word $U$ over  the positive alphabet $\{ a, b, c,
d\}$ such that $|U| \le |W| \le x$ and any cyclic permutation
of $U$ contains no subwords of the form $a^2$ and $g_1g_2$,
where $g_1, g_2 \in \{ b, c, d \}$. Indeed, at most $x$ copies
of relators $a^{2}$, $b^{2}$, $c^{2}$, $d^{2}$ are sufficient
to turn $W$ into a positive word $W'$ with $|W'| \le |W|$.
Then, decreasing $|W'|$, we can delete the unwanted subwords by
applying $\le 2x$ additional copies of relators $a^{2}$,
$b^{2}$, $c^{ 2}$, $d^{2}$, $bcd^{}$.

Clearly, $U \in \N(\infty)$. Let $\psi_0(U) = (U_0, U_1)$. By
the definitions of $\psi_0$ and $U$, we have $\max(|U_0|,
|U_1|) \le |U|/2 \le x/2$.  It follows from Lemmas~3.1--3.2
that $U_0, U_1 \in \N(\infty)$ and, in the free group
$F(a,b,c,d)$,
 \begin{equation}\label{UaV}
U = a \sigma(U_0) a \sigma(U_1) V \ ,
 \end{equation}
where $V \in \N(2)$. Since $|\sigma(U_0)|  \le 3 |U_0|$,
$|\sigma(U_1)|  \le 3 |U_1|$, we can assume that
   \begin{equation*}
 |V| \le 2\cdot 3|U_0| + |U| +2 \le 4x +2 \ .
   \end{equation*}

It is known \cite{GH}  that the group $\Ga(2) =
F(a,b,c,d)/\N(2)$ contains a subgroup of finite index which is
isomorphic to the direct product $F_2\times F_2$, where $F_2$
is a free group of rank 2. Moreover,  $\Ga(i) =
F(a,b,c,d)/\N(i)$, where $i \ge 2$ is finite, contains a
subgroup of finite index isomorphic to the direct product of
$2^i$ copies of $F_2$, see \cite{GH}. Since the equivalence
class of Dehn $2$-functions (for finite presentations) does not
change when passing to subgroups of finite index, and the Dehn
$2$-function of $F_2\times F_2$ is quadratic,  it follows that
there exists a constant $C_0>0$ for the presentation $\Ga(2)$
such that if $Y \in \N(2)$ then $Y$ can be written as a product
of at most $C_0 |Y|^2$  conjugates of elements of $\R(2)^{ \pm
1}$. In particular, the word  $V$ can be written as a product
of at most
 \begin{equation*}
 C_0 |V|^2  \le  C_0 (4x +2)^2
 \end{equation*}
 conjugates of elements of $\R(2)^{ \pm 1}$. Since $\R(2)
 \subset \s^*(0) \bigcup \s^*(1)$, it follows that $V$ possesses an
 $h^*$-tuple  $\tau^V = (\tau^V_0,
\tau^V_1, \dots )$ such that
\begin{equation}\label{V}
\tau^V_0+ \tau^V_1  \le  C_0 (4x +2)^2 \  \quad \text{and}
\quad \tau^V_i = 0 \quad \text{for} \quad   i>1 \ .
 \end{equation}

By the induction hypothesis applied to the words $U_0$, $U_1$,
we obtain the existence of $h^*$-tuples  $\tau^k = (\tau^k_0,
\tau^k_1, \dots )$ for $U_k$, where $k =0,1$, such that
$\tau^k_i = 0$ if $i > \log_2(x/2)$ and
\begin{equation}\label{IH}
\tau^k_0 \le C(x/2)^2 \  , \qquad \tau_i^k \le \frac
{C}{2^{i-1}} (x/2)^2
\end{equation}
for $i$ with $1 \le i \le \log_2(x/2)$.

Note that if
\begin{equation*}
U_k =  \prod_{j=1}^{L_k}  X_{jk}  R_{jk}^{\e_{jk}}  X_{jk}^{-1}
\ ,
\end{equation*}
where $X_{jk} \in F(a,b,c,d)$, $R_{jk} \in \R^*(\infty)$,
$\e_{jk} = \pm 1$, $k= 0,1$,  then
\begin{equation*}
\sigma(U_k) =  \prod_{j=1}^{L_k}  \sigma(X_{jk})
\sigma(R_{jk})^{\e_{jk}} \sigma(X_{jk})^{-1} \ ,
\end{equation*}
where $\sigma(X_{jk}) \in F(a,b,c,d)$. In addition, up to a
cyclic permutation (in case $R_{jk} = bcd$), we have that
$\sigma(R_{jk}) \in \R^*(\infty)$ with either
$h(\sigma(R_{jk})) = h(R_{jk})$, if $R_{jk} \in \{b^{2}, c^{
2}, d^{2}, bcd \}$, or $h(\sigma(R_{jk})) = h(R_{jk}) +1$
otherwise. As was pointed out above, $\sigma(a^2) \not\in
\R(\infty)$ and this is the reason for extending the set
$\R(\infty)$ to $\R^*(\infty)$.

Referring to the presentation \eqref{UaV} for $U$, we observe
that $U$ is a product of words $a^2$, $a^{-1} \sigma(U_0) a$,
$\sigma(U_1)$, $V$. Hence, an  $h^*$-tuple  for $U$ can be
obtained as the sum of
 $h^*$-tuples for $a^2$, $a^{-1} \sigma(U_0) a$,
$\sigma(U_1)$, $V$.  This observation implies the existence of
an $h^*$-tuple $\tau^U = (\tau^U_0, \tau^U_1, \dots )$ for $U$
with the following properties.

First, since $\tau^k_{i'} = 0$, where $i' > \log_2(x/2)$, $k =0,1$, and
$\tau^V_{i'} = 0$, where $i' > 1$,  it follows that $\tau_{i}^U \le
\tau^0_{i-1}+\tau^1_{i-1} + \tau^V_{i} =0$ when $i > \log_2x =
\log_2(x/2) +1 \ge 1$ for $x \ge 2$. Second, in view of the estimates
\eqref{V}, \eqref{IH}, $x \ge 2$ and the definitions, we have
\begin{align*}
 \tau_0^U & \le 1 + \tau^V_{0}+  \tau^0_0 + \tau^1_0 \le C_0(4x+2)^2
+ 1+ 2C(x/2)^2 \\ & \le   C_0(4x+2)^2 + \frac C2x^2 + 1 \le
C_0(5x)^2 +1+ \frac C2 x^2 \ , \\
\tau_1^U & \le \tau^V_{1}   + \tau^0_0 + \tau^1_0 \le
C_0(4x+2)^2 + 2C(x/2)^2
 \le   C_0(5x)^2 +
\frac C2x^2  \ , \quad \text{and} \\
\tau_i^U & \le \tau^0_{i-1} + \tau^1_{i-1} \le 2 \cdot
\frac{C}{2^{i-2}}(x/2)^2 = \frac{C}{2^{i-1}} x^2 \ , \qquad
\quad \text{where} \quad i \ge 2 \ .
\end{align*}
Choosing $C \ge 50(C_0+1)$, we obtain that $\tau_0^U$,
$\tau_1^U < Cx^2 - 3x$. Recall that $U$ was obtained from $W$
by at most $3x$ applications of copies of relators $a^{2}$,
$b^{2}$, $c^{2}$, $d^{2}$, $bcd^{}$.  This remark finally
proves the existence of an $h^*$-tuple  $\tau^W =(\tau_0^W,
\tau_1^W, \dots)$ for $W$ with all of the desired properties.

The inequality  \eqref{2inq1} is now obvious from the proven
estimates for $\tau_0^W, \tau_1^W,  \dots$. Note that if $R \in
\s^*(i)$ then
\begin{equation}\label{RE}
|R| \le | \sigma^i( (adacac)^4 ) | = 12 \cdot 2^{i +1} = 3
\cdot 2^{i+3} \ ,
\end{equation}
because every application of $\sigma$ doubles the length. The
 inequality \eqref{2inq2} also becomes evident.
\end{proof}

{\em Proof of Theorem~1.8:}  Let $W \in \N(\infty)$ and $|W| \le x$.
According to Lemma~3.3, there exists an $h^*$-tuple $\tau^W =(\tau_0^W,
\tau_1^W, \dots)$ for $W$ with the described
 properties.  Since an $h^*$-tuple  is defined  with respect to
the relator set $\R^*(\infty)$, we may have relators
$\sigma^i(a^2)$ in  a presentation for $W$ of the form
\eqref{DF2}, $1 \le i \le \log_2x$, which are not in
$\R(\infty)$. Recall that $\sigma^i(a^2)$ belongs to the normal
closure of words $a^2, b^2$, $c^2, d^2$. We also note that, as
follows from the definition of $\sigma$, see \eqref{sig},
$\sigma^i(a^2)$ is a positive word and has length $2^{i+2}-2$.

Therefore, we can turn $\sigma^i(a^2)$ into the empty word by
consecutive deletions of $\frac{2^{i+2}-2}{2} =2^{i+1}-1$
subwords  of the form $g^2$, $g \in \{a,b,c,d\}$. This means
that $\sigma^i(a^2)$ is a product of $2^{i+1}-1$ conjugates of
relators $a^2, b^2$, $c^2, d^2$.

Consequently,  replacing each relator $\sigma^i(a^2)$ (or its
inverse), $i \ge 1$, by a product of $2^{i+1}-1$ conjugates of
relators $a^2$, $b^2$, $c^2$, $d^2$ (or their inverses), we can
turn a product \eqref{DF2}, defined for $W$ by means of the
extended set $\R^*(\infty)$, into a similar product,  defined
for $W$ by means of the original set $\R(\infty)$. Hence, there
exists an $h$-tuple $ \tau' = (\tau'_0, \tau'_1, \dots)$,
defined for $W$ by means of $\R(\infty)$, such that
\begin{align}\notag
\tau'_0 & =  \tau_0^W + \sum_{i=1}^{[ \log_2x ]} (2^{i+1}-1)\tau_i^W  \le
Cx^2 +  \sum_{i=1}^{[ \log_2x ]} (2^{i+1}-1) \cdot \frac{C} {2^{i-1}} x^2
\\ & \le  Cx^2 + \sum_{i=1}^{[ \log_2x ]}
(4- 2^{1-i}) C x^2 \le 4 C  x^2 \log_2x \ ,  \label{tau1} \\
\tau'_i & \le \tau_i^W \le \frac {C}{2^{i-1}} x^2   \quad
\text{if} \quad 1 \le i \le  \log_2 x \ ,  \quad \text{and}
\quad \tau'_i  = 0 \ \quad \text{if}  \quad i > \log_2 x \ ,
\label{tau2}
\end{align}
where $[ \log_2x ]$ is the greatest integer $\ell$ with $\ell
\le \log_2x$.

Let $\D$ be a van Kampen diagram over the presentation
 \eqref{pr1}   with $|\p \D | \le x$. Assuming that $\D$ is
minimal relative to $|\D(2)|$, in view of estimates \eqref{RE},
\eqref{tau1}, \eqref{tau2} and Lemma~3.3, we can estimate the
numbers $|\D(j)|$ of $j$-cells in $\D$, $j =1,2$, and, hence,
the Dehn functions $f_{j,\Ga}(x)$  of the presentation
\eqref{pr1}  as follows
\begin{align*}
f_{1,\Ga}(x) & \le 2\tau'_0 +
 \sum_{i=0}^{[ \log_2x ]}3 \cdot 2^{i+3} \tau'_i + x/2
 \le 8C x^2 \log_2x + 3\cdot 2^4C x^2 \log_2x + x/2 \\ &
 \le 2^6 C x^2 \log_2x + x/2 \ , \\
f_{2,\Ga}(x) & \le \sum_{i=0}^{[ \log_2x ]} \tau'_i \le 4C x^2
\log_2x  + \sum_{i=1}^{[ \log_2x ]} \tau_i
   \le 4 C  x^2 \log_2x + 2Cx^2  \ .
\end{align*}
Since $f_{1, \Ga}(1) = f_{2,\Ga}(1) =0$ and $C>1$, it follows from these
estimates that
\begin{equation*}
f_{1,\Ga}(x)  \le 55 C x^2 \log_2x \quad  \mbox{and} \quad  f_{2,\Ga}(x)
\le 6 C x^2 \log_2x \ .
\end{equation*}
By Theorem~1.4(a), $f_{0,\Ga}(x) \le 2 f_{1,\Ga}(x)$, hence we can pick
$C_1 = 110 C$ and the proof of Theorem~1.8 is complete. \qed
\medskip

{\em Proof of Corollary~1.9:} Let $W$ be a word in the normal closure
$\langle \langle \R_t \rangle \rangle$  of the relator set $\R_t$ of the
presentation  \eqref{pr1a} of $\Ga_t$. By $|W|_\ell$ denote the number of
occurrences of letters $\ell, \ell^{-1}$ in $W$, where $\ell \in \{
a,b,c,d,t \}$. Since $t g t^{-1} = \sigma(g)$, where $g \in \{ a,b,c,d
\}$, it follows from the definitions, see \eqref{sig}, \eqref{pr1a}, that
we can apply relations $t g t^{-1} = \sigma(g)$ to the cyclic word $W$
and eliminate all occurrences of $t, t^{-1}$. Note that each application
of $t a t^{-1} = a ca$  adds two extra letters and doubles the number of
occurrences of $a$, whereas an application of $t g t^{-1} = \sigma(g)$,
where $g \in \{b,c,d \}$, preserves the length. Hence, when eliminating
two successive  occurrences of $t$, $t^{-1}$ in the cyclic word $W$, we
get a new cyclic word $W_1$ such that
\begin{align*}
|W_1|_t = |W|_t - 2 \ , \quad |W_1|_a \le 2 |W|_a  \ , \quad
\text{and} \quad |W_1| \le  |W|+ 2|W|_a \ .
\end{align*}

Iterating, in $k\ge 1$ steps, we get a word $W_k$ such that
\begin{align}\notag
|W_k|_t & = |W|_t - 2k \ , \quad  |W_k|_a \le 2^k |W|_a  \ ,
\quad \text{and} \quad \\ |W_k|  & \le  |W|+ |W|_a(2^1 + \dots
+ 2^k) =  |W|+ |W|_a(2^{k+1} -2) \ . \label{dW}
\end{align}

Making $|W|_t/2$ such steps, we get a word $U = W_{|W|_t/2 }$
with no occurrences of $t, t^{-1}$. Recall that $W \in \langle
\langle \R_t \rangle \rangle$, whence $|W|_t$ is even.

Note that, in view of inequalities $|W|_a +|W|_t \le |W| $ and
$y < 2^y$, we have
\begin{align}\label{Wa}
 |W|_a \cdot 2^{|W|_t/2}   \le    |W|_a \cdot 2^{(|W| - |W|_a )/2}
 \le \frac{|W|_a}{2^{|W|_a/2} }\cdot 2^{|W|/2} < 2\cdot 2^{|W|/2} \
 .
\end{align}

Hence, according to \eqref{dW}, \eqref{Wa},
\begin{align}\notag
|U|  & \le  |W|+ |W|_a\cdot (2^{|W|_t/2 +1} -2)  \le    |W|+
2|W|_a\cdot 2^{|W|_t/2}  \\    & \le    2\cdot 2^{|W| /2} +
4\cdot 2^{|W| /2}  =  6\cdot 2^{|W| /2}  \ . \label{dU}
\end{align}

It follows from inequalities \eqref{dW}, \eqref{Wa} and $y^2 <
2^3 \cdot 2^{y/2}$, where $y > 0$, that the number of relations
of the form $t g t^{-1} = \sigma(g)$, where $g \in \{
a,b,c,d\}$,  that are needed to obtain $U$ from $W$, does not
exceed the following sum
\begin{align}\notag
|W|+ |W_1| + \dots + & |W_{|W|_t/2 }|  \le \tfrac 12 |W|^2 +
|W|_a \cdot \sum_{k=1}^{|W|_t/2} 2^{k+1}  \\
& \le  2^2\cdot 2^{|W|/2}  +|W|_a \cdot 2^{|W|_t/2 +2}  \le
12\cdot 2^{|W|/2} \ . \label{dU2}
\end{align}

As we saw in  the proof of Theorem~1.8,  the word $U$ has an $h$-tuple
$\tau' = (\tau'_0, \tau'_1, \dots)$ defined for $U$ by means of
$\R(\infty)$ whose entries satisfy the inequalities \eqref{tau1},
\eqref{tau2} in which $x = |U|$.

 Note that if $R = \sigma^j(V)$, where $V \in \{
(ad)^4, (adacac)^4\}$,  $j \ge 1$, then there exists a van
Kampen diagram $\E_R$ over the presentation \eqref{pr1a} of
$\Ga_t$ such that  $\ph(\p \E_R) \equiv R$, $\E_R$ contains a
face $\Pi$, with  $\ph( \p \Pi) \equiv V^{-1}$, and $\Pi$ is
surrounded by $j$ annuli each of which consists of faces
corresponding to relators $t g t^{-1} \sigma(g)^{-1}$, $g \in
\{a,b,c,d\}$, see Figure~3.

\unitlength 1mm \linethickness{0.8pt}
\ifx\plotpoint\undefined\newsavebox{\plotpoint}\fi 
\begin{picture}(70.00,47.75)(-10.00,5.00)
\put(48.00,32.50){\circle{10.00}}
\put(58.00,32.50){\line(0,1){0.54}}
\put(57.99,33.04){\line(0,1){0.54}}
\put(57.94,33.58){\line(0,1){0.54}}
\put(57.87,34.12){\line(0,1){0.53}}
\multiput(57.77,34.65)(-0.07,0.26){2}{\line(0,1){0.26}}
\multiput(57.64,35.18)(-0.08,0.26){2}{\line(0,1){0.26}}
\multiput(57.48,35.69)(-0.09,0.25){2}{\line(0,1){0.25}}
\multiput(57.29,36.20)(-0.11,0.25){2}{\line(0,1){0.25}}
\multiput(57.08,36.70)(-0.08,0.16){3}{\line(0,1){0.16}}
\multiput(56.84,37.18)(-0.09,0.16){3}{\line(0,1){0.16}}
\multiput(56.57,37.66)(-0.10,0.15){3}{\line(0,1){0.15}}
\multiput(56.28,38.11)(-0.11,0.15){3}{\line(0,1){0.15}}
\multiput(55.96,38.55)(-0.11,0.14){3}{\line(0,1){0.14}}
\multiput(55.62,38.97)(-0.09,0.10){4}{\line(0,1){0.10}}
\multiput(55.26,39.38)(-0.10,0.10){4}{\line(0,1){0.10}}
\multiput(54.88,39.76)(-0.10,0.09){4}{\line(-1,0){0.10}}
\multiput(54.47,40.12)(-0.14,0.11){3}{\line(-1,0){0.14}}
\multiput(54.05,40.46)(-0.15,0.11){3}{\line(-1,0){0.15}}
\multiput(53.61,40.78)(-0.15,0.10){3}{\line(-1,0){0.15}}
\multiput(53.16,41.07)(-0.16,0.09){3}{\line(-1,0){0.16}}
\multiput(52.68,41.34)(-0.16,0.08){3}{\line(-1,0){0.16}}
\multiput(52.20,41.58)(-0.25,0.11){2}{\line(-1,0){0.25}}
\multiput(51.70,41.79)(-0.25,0.09){2}{\line(-1,0){0.25}}
\multiput(51.19,41.98)(-0.26,0.08){2}{\line(-1,0){0.26}}
\multiput(50.68,42.14)(-0.26,0.07){2}{\line(-1,0){0.26}}
\put(50.15,42.27){\line(-1,0){0.53}}
\put(49.62,42.37){\line(-1,0){0.54}}
\put(49.08,42.44){\line(-1,0){0.54}}
\put(48.54,42.49){\line(-1,0){1.08}}
\put(47.46,42.49){\line(-1,0){0.54}}
\put(46.92,42.44){\line(-1,0){0.54}}
\put(46.38,42.37){\line(-1,0){0.53}}
\multiput(45.85,42.27)(-0.26,-0.07){2}{\line(-1,0){0.26}}
\multiput(45.32,42.14)(-0.26,-0.08){2}{\line(-1,0){0.26}}
\multiput(44.81,41.98)(-0.25,-0.09){2}{\line(-1,0){0.25}}
\multiput(44.30,41.79)(-0.25,-0.11){2}{\line(-1,0){0.25}}
\multiput(43.80,41.58)(-0.16,-0.08){3}{\line(-1,0){0.16}}
\multiput(43.32,41.34)(-0.16,-0.09){3}{\line(-1,0){0.16}}
\multiput(42.84,41.07)(-0.15,-0.10){3}{\line(-1,0){0.15}}
\multiput(42.39,40.78)(-0.15,-0.11){3}{\line(-1,0){0.15}}
\multiput(41.95,40.46)(-0.14,-0.11){3}{\line(-1,0){0.14}}
\multiput(41.53,40.12)(-0.10,-0.09){4}{\line(-1,0){0.10}}
\multiput(41.12,39.76)(-0.10,-0.10){4}{\line(0,-1){0.10}}
\multiput(40.74,39.38)(-0.09,-0.10){4}{\line(0,-1){0.10}}
\multiput(40.38,38.97)(-0.11,-0.14){3}{\line(0,-1){0.14}}
\multiput(40.04,38.55)(-0.11,-0.15){3}{\line(0,-1){0.15}}
\multiput(39.72,38.11)(-0.10,-0.15){3}{\line(0,-1){0.15}}
\multiput(39.43,37.66)(-0.09,-0.16){3}{\line(0,-1){0.16}}
\multiput(39.16,37.18)(-0.08,-0.16){3}{\line(0,-1){0.16}}
\multiput(38.92,36.70)(-0.11,-0.25){2}{\line(0,-1){0.25}}
\multiput(38.71,36.20)(-0.09,-0.25){2}{\line(0,-1){0.25}}
\multiput(38.52,35.69)(-0.08,-0.26){2}{\line(0,-1){0.26}}
\multiput(38.36,35.18)(-0.07,-0.26){2}{\line(0,-1){0.26}}
\put(38.23,34.65){\line(0,-1){0.53}}
\put(38.13,34.12){\line(0,-1){0.54}}
\put(38.06,33.58){\line(0,-1){0.54}}
\put(38.01,33.04){\line(0,-1){1.62}}
\put(38.06,31.42){\line(0,-1){0.54}}
\put(38.13,30.88){\line(0,-1){0.53}}
\multiput(38.23,30.35)(0.07,-0.26){2}{\line(0,-1){0.26}}
\multiput(38.36,29.82)(0.08,-0.26){2}{\line(0,-1){0.26}}
\multiput(38.52,29.31)(0.09,-0.25){2}{\line(0,-1){0.25}}
\multiput(38.71,28.80)(0.11,-0.25){2}{\line(0,-1){0.25}}
\multiput(38.92,28.30)(0.08,-0.16){3}{\line(0,-1){0.16}}
\multiput(39.16,27.82)(0.09,-0.16){3}{\line(0,-1){0.16}}
\multiput(39.43,27.34)(0.10,-0.15){3}{\line(0,-1){0.15}}
\multiput(39.72,26.89)(0.11,-0.15){3}{\line(0,-1){0.15}}
\multiput(40.04,26.45)(0.11,-0.14){3}{\line(0,-1){0.14}}
\multiput(40.38,26.03)(0.09,-0.10){4}{\line(0,-1){0.10}}
\multiput(40.74,25.62)(0.10,-0.10){4}{\line(0,-1){0.10}}
\multiput(41.12,25.24)(0.10,-0.09){4}{\line(1,0){0.10}}
\multiput(41.53,24.88)(0.14,-0.11){3}{\line(1,0){0.14}}
\multiput(41.95,24.54)(0.15,-0.11){3}{\line(1,0){0.15}}
\multiput(42.39,24.22)(0.15,-0.10){3}{\line(1,0){0.15}}
\multiput(42.84,23.93)(0.16,-0.09){3}{\line(1,0){0.16}}
\multiput(43.32,23.66)(0.16,-0.08){3}{\line(1,0){0.16}}
\multiput(43.80,23.42)(0.25,-0.11){2}{\line(1,0){0.25}}
\multiput(44.30,23.21)(0.25,-0.09){2}{\line(1,0){0.25}}
\multiput(44.81,23.02)(0.26,-0.08){2}{\line(1,0){0.26}}
\multiput(45.32,22.86)(0.26,-0.07){2}{\line(1,0){0.26}}
\put(45.85,22.73){\line(1,0){0.53}}
\put(46.38,22.63){\line(1,0){0.54}}
\put(46.92,22.56){\line(1,0){0.54}}
\put(47.46,22.51){\line(1,0){1.08}}
\put(48.54,22.51){\line(1,0){0.54}}
\put(49.08,22.56){\line(1,0){0.54}}
\put(49.62,22.63){\line(1,0){0.53}}
\multiput(50.15,22.73)(0.26,0.07){2}{\line(1,0){0.26}}
\multiput(50.68,22.86)(0.26,0.08){2}{\line(1,0){0.26}}
\multiput(51.19,23.02)(0.25,0.09){2}{\line(1,0){0.25}}
\multiput(51.70,23.21)(0.25,0.11){2}{\line(1,0){0.25}}
\multiput(52.20,23.42)(0.16,0.08){3}{\line(1,0){0.16}}
\multiput(52.68,23.66)(0.16,0.09){3}{\line(1,0){0.16}}
\multiput(53.16,23.93)(0.15,0.10){3}{\line(1,0){0.15}}
\multiput(53.61,24.22)(0.15,0.11){3}{\line(1,0){0.15}}
\multiput(54.05,24.54)(0.14,0.11){3}{\line(1,0){0.14}}
\multiput(54.47,24.88)(0.10,0.09){4}{\line(1,0){0.10}}
\multiput(54.88,25.24)(0.10,0.10){4}{\line(0,1){0.10}}
\multiput(55.26,25.62)(0.09,0.10){4}{\line(0,1){0.10}}
\multiput(55.62,26.03)(0.11,0.14){3}{\line(0,1){0.14}}
\multiput(55.96,26.45)(0.11,0.15){3}{\line(0,1){0.15}}
\multiput(56.28,26.89)(0.10,0.15){3}{\line(0,1){0.15}}
\multiput(56.57,27.34)(0.09,0.16){3}{\line(0,1){0.16}}
\multiput(56.84,27.82)(0.08,0.16){3}{\line(0,1){0.16}}
\multiput(57.08,28.30)(0.11,0.25){2}{\line(0,1){0.25}}
\multiput(57.29,28.80)(0.09,0.25){2}{\line(0,1){0.25}}
\multiput(57.48,29.31)(0.08,0.26){2}{\line(0,1){0.26}}
\multiput(57.64,29.82)(0.07,0.26){2}{\line(0,1){0.26}}
\put(57.77,30.35){\line(0,1){0.53}}
\put(57.87,30.88){\line(0,1){0.54}}
\put(57.94,31.42){\line(0,1){1.08}}
\put(63.00,32.50){\line(0,1){0.75}}
\put(62.98,33.25){\line(0,1){0.75}}
\put(62.93,33.99){\line(0,1){0.74}}
\multiput(62.83,34.74)(-0.06,0.37){2}{\line(0,1){0.37}}
\multiput(62.70,35.47)(-0.08,0.36){2}{\line(0,1){0.36}}
\multiput(62.54,36.20)(-0.10,0.36){2}{\line(0,1){0.36}}
\multiput(62.33,36.92)(-0.12,0.35){2}{\line(0,1){0.35}}
\multiput(62.10,37.63)(-0.09,0.23){3}{\line(0,1){0.23}}
\multiput(61.82,38.33)(-0.10,0.23){3}{\line(0,1){0.23}}
\multiput(61.51,39.01)(-0.11,0.22){3}{\line(0,1){0.22}}
\multiput(61.17,39.67)(-0.09,0.16){4}{\line(0,1){0.16}}
\multiput(60.80,40.32)(-0.10,0.16){4}{\line(0,1){0.16}}
\multiput(60.39,40.95)(-0.11,0.15){4}{\line(0,1){0.15}}
\multiput(59.96,41.56)(-0.12,0.15){4}{\line(0,1){0.15}}
\multiput(59.49,42.14)(-0.10,0.11){5}{\line(0,1){0.11}}
\multiput(59.00,42.70)(-0.10,0.11){5}{\line(0,1){0.11}}
\multiput(58.47,43.24)(-0.11,0.10){5}{\line(-1,0){0.11}}
\multiput(57.93,43.75)(-0.11,0.10){5}{\line(-1,0){0.11}}
\multiput(57.35,44.23)(-0.15,0.11){4}{\line(-1,0){0.15}}
\multiput(56.76,44.68)(-0.15,0.11){4}{\line(-1,0){0.15}}
\multiput(56.14,45.10)(-0.16,0.10){4}{\line(-1,0){0.16}}
\multiput(55.50,45.49)(-0.22,0.12){3}{\line(-1,0){0.22}}
\multiput(54.84,45.85)(-0.22,0.11){3}{\line(-1,0){0.22}}
\multiput(54.17,46.17)(-0.23,0.10){3}{\line(-1,0){0.23}}
\multiput(53.48,46.46)(-0.23,0.09){3}{\line(-1,0){0.23}}
\multiput(52.78,46.72)(-0.36,0.11){2}{\line(-1,0){0.36}}
\multiput(52.06,46.94)(-0.36,0.09){2}{\line(-1,0){0.36}}
\multiput(51.34,47.12)(-0.37,0.07){2}{\line(-1,0){0.37}}
\put(50.60,47.27){\line(-1,0){0.74}}
\put(49.87,47.38){\line(-1,0){0.74}}
\put(49.12,47.46){\line(-1,0){0.75}}
\put(48.37,47.50){\line(-1,0){0.75}}
\put(47.63,47.50){\line(-1,0){0.75}}
\put(46.88,47.46){\line(-1,0){0.74}}
\put(46.13,47.38){\line(-1,0){0.74}}
\multiput(45.40,47.27)(-0.37,-0.07){2}{\line(-1,0){0.37}}
\multiput(44.66,47.12)(-0.36,-0.09){2}{\line(-1,0){0.36}}
\multiput(43.94,46.94)(-0.36,-0.11){2}{\line(-1,0){0.36}}
\multiput(43.22,46.72)(-0.23,-0.09){3}{\line(-1,0){0.23}}
\multiput(42.52,46.46)(-0.23,-0.10){3}{\line(-1,0){0.23}}
\multiput(41.83,46.17)(-0.22,-0.11){3}{\line(-1,0){0.22}}
\multiput(41.16,45.85)(-0.22,-0.12){3}{\line(-1,0){0.22}}
\multiput(40.50,45.49)(-0.16,-0.10){4}{\line(-1,0){0.16}}
\multiput(39.86,45.10)(-0.15,-0.11){4}{\line(-1,0){0.15}}
\multiput(39.24,44.68)(-0.15,-0.11){4}{\line(-1,0){0.15}}
\multiput(38.65,44.23)(-0.11,-0.10){5}{\line(-1,0){0.11}}
\multiput(38.07,43.75)(-0.11,-0.10){5}{\line(-1,0){0.11}}
\multiput(37.53,43.24)(-0.10,-0.11){5}{\line(0,-1){0.11}}
\multiput(37.00,42.70)(-0.10,-0.11){5}{\line(0,-1){0.11}}
\multiput(36.51,42.14)(-0.12,-0.15){4}{\line(0,-1){0.15}}
\multiput(36.04,41.56)(-0.11,-0.15){4}{\line(0,-1){0.15}}
\multiput(35.61,40.95)(-0.10,-0.16){4}{\line(0,-1){0.16}}
\multiput(35.20,40.32)(-0.09,-0.16){4}{\line(0,-1){0.16}}
\multiput(34.83,39.67)(-0.11,-0.22){3}{\line(0,-1){0.22}}
\multiput(34.49,39.01)(-0.10,-0.23){3}{\line(0,-1){0.23}}
\multiput(34.18,38.33)(-0.09,-0.23){3}{\line(0,-1){0.23}}
\multiput(33.90,37.63)(-0.12,-0.35){2}{\line(0,-1){0.35}}
\multiput(33.67,36.92)(-0.10,-0.36){2}{\line(0,-1){0.36}}
\multiput(33.46,36.20)(-0.08,-0.36){2}{\line(0,-1){0.36}}
\multiput(33.30,35.47)(-0.06,-0.37){2}{\line(0,-1){0.37}}
\put(33.17,34.74){\line(0,-1){0.74}}
\put(33.07,33.99){\line(0,-1){0.75}}
\put(33.02,33.25){\line(0,-1){2.24}}
\put(33.07,31.01){\line(0,-1){0.74}}
\multiput(33.17,30.26)(0.06,-0.37){2}{\line(0,-1){0.37}}
\multiput(33.30,29.53)(0.08,-0.36){2}{\line(0,-1){0.36}}
\multiput(33.46,28.80)(0.10,-0.36){2}{\line(0,-1){0.36}}
\multiput(33.67,28.08)(0.12,-0.35){2}{\line(0,-1){0.35}}
\multiput(33.90,27.37)(0.09,-0.23){3}{\line(0,-1){0.23}}
\multiput(34.18,26.67)(0.10,-0.23){3}{\line(0,-1){0.23}}
\multiput(34.49,25.99)(0.11,-0.22){3}{\line(0,-1){0.22}}
\multiput(34.83,25.33)(0.09,-0.16){4}{\line(0,-1){0.16}}
\multiput(35.20,24.68)(0.10,-0.16){4}{\line(0,-1){0.16}}
\multiput(35.61,24.05)(0.11,-0.15){4}{\line(0,-1){0.15}}
\multiput(36.04,23.44)(0.12,-0.15){4}{\line(0,-1){0.15}}
\multiput(36.51,22.86)(0.10,-0.11){5}{\line(0,-1){0.11}}
\multiput(37.00,22.30)(0.10,-0.11){5}{\line(0,-1){0.11}}
\multiput(37.53,21.76)(0.11,-0.10){5}{\line(1,0){0.11}}
\multiput(38.07,21.25)(0.11,-0.10){5}{\line(1,0){0.11}}
\multiput(38.65,20.77)(0.15,-0.11){4}{\line(1,0){0.15}}
\multiput(39.24,20.32)(0.15,-0.11){4}{\line(1,0){0.15}}
\multiput(39.86,19.90)(0.16,-0.10){4}{\line(1,0){0.16}}
\multiput(40.50,19.51)(0.22,-0.12){3}{\line(1,0){0.22}}
\multiput(41.16,19.15)(0.22,-0.11){3}{\line(1,0){0.22}}
\multiput(41.83,18.83)(0.23,-0.10){3}{\line(1,0){0.23}}
\multiput(42.52,18.54)(0.23,-0.09){3}{\line(1,0){0.23}}
\multiput(43.22,18.28)(0.36,-0.11){2}{\line(1,0){0.36}}
\multiput(43.94,18.06)(0.36,-0.09){2}{\line(1,0){0.36}}
\multiput(44.66,17.88)(0.37,-0.07){2}{\line(1,0){0.37}}
\put(45.40,17.73){\line(1,0){0.74}}
\put(46.13,17.62){\line(1,0){0.74}}
\put(46.88,17.54){\line(1,0){0.75}}
\put(47.63,17.50){\line(1,0){0.75}}
\put(48.37,17.50){\line(1,0){0.75}}
\put(49.12,17.54){\line(1,0){0.74}}
\put(49.87,17.62){\line(1,0){0.74}}
\multiput(50.60,17.73)(0.37,0.07){2}{\line(1,0){0.37}}
\multiput(51.34,17.88)(0.36,0.09){2}{\line(1,0){0.36}}
\multiput(52.06,18.06)(0.36,0.11){2}{\line(1,0){0.36}}
\multiput(52.78,18.28)(0.23,0.09){3}{\line(1,0){0.23}}
\multiput(53.48,18.54)(0.23,0.10){3}{\line(1,0){0.23}}
\multiput(54.17,18.83)(0.22,0.11){3}{\line(1,0){0.22}}
\multiput(54.84,19.15)(0.22,0.12){3}{\line(1,0){0.22}}
\multiput(55.50,19.51)(0.16,0.10){4}{\line(1,0){0.16}}
\multiput(56.14,19.90)(0.15,0.11){4}{\line(1,0){0.15}}
\multiput(56.76,20.32)(0.15,0.11){4}{\line(1,0){0.15}}
\multiput(57.35,20.77)(0.11,0.10){5}{\line(1,0){0.11}}
\multiput(57.93,21.25)(0.11,0.10){5}{\line(1,0){0.11}}
\multiput(58.47,21.76)(0.10,0.11){5}{\line(0,1){0.11}}
\multiput(59.00,22.30)(0.10,0.11){5}{\line(0,1){0.11}}
\multiput(59.49,22.86)(0.12,0.15){4}{\line(0,1){0.15}}
\multiput(59.96,23.44)(0.11,0.15){4}{\line(0,1){0.15}}
\multiput(60.39,24.05)(0.10,0.16){4}{\line(0,1){0.16}}
\multiput(60.80,24.68)(0.09,0.16){4}{\line(0,1){0.16}}
\multiput(61.17,25.33)(0.11,0.22){3}{\line(0,1){0.22}}
\multiput(61.51,25.99)(0.10,0.23){3}{\line(0,1){0.23}}
\multiput(61.82,26.67)(0.09,0.23){3}{\line(0,1){0.23}}
\multiput(62.10,27.37)(0.12,0.35){2}{\line(0,1){0.35}}
\multiput(62.33,28.08)(0.10,0.36){2}{\line(0,1){0.36}}
\multiput(62.54,28.80)(0.08,0.36){2}{\line(0,1){0.36}}
\multiput(62.70,29.53)(0.06,0.37){2}{\line(0,1){0.37}}
\put(62.83,30.26){\line(0,1){0.74}}
\put(62.93,31.01){\line(0,1){1.49}}
\put(48.00,37.50){\line(0,1){10.00}}
\put(48.00,27.50){\line(0,-1){9.75}}
\put(53.00,32.50){\line(1,0){10.00}}
\put(43.,32.50){\line(-1,0){9.9}}
\put(40.75,39.25){\line(-1,1){3.750}}
\put(55.00,25.50){\line(1,-1){3.75}}
\put(55.00,39.75){\line(1,1){3.50}}
\put(40.75,25.25){\line(-1,-1){3.50}}
\put(48.00,46.75){\vector(0,-1){3.00}}
\put(48.00,41.50){\vector(0,-1){2.75}}
\put(38.25,41.75){\vector(1,-1){1.75}}
\put(57.75,42.25){\vector(-1,-1){2.25}}
\put(34.00,32.50){\vector(1,0){2.50}}
\put(39.50,32.50){\vector(1,0){2.25}}
\put(62.00,32.50){\vector(-1,0){2.50}}
\put(56.50,32.50){\vector(-1,0){2.25}}
\put(48.00,18.75){\vector(0,1){2.25}}
\put(48.00,24.00){\vector(0,1){2.25}}
\put(38.50,23.00){\vector(1,1){1.50}}
\put(57.75,22.75){\vector(-1,1){1.75}}
\put(46.800,31.50){$\Pi$}
 \put(49.5,44.00){$t$} \put(57.75,39.250){$t$}
 \put(59.75,29.25){$t$} \put(54.75,22.00){$t$}
 \put(45.5,19.200){$t$} \put(37.25,25.00){$t$}
 \put(35.00,34.){$t$} \put(40.100,41.5){$t$}
 \put(49.75,39.0){$t$} \put(54.5,29.250){$t$}
 \put(45.500,24.5){$t$} \put(40.50,34.00){$t$}
 \put(43.00,10.00){Figure 3} \put(68.00,30.00){$\E_R$}
 \put(58.75,43.25){\line(-1,-1){3.75}}
 \put(37.25,21.75){\line(1,1){3.75}}
\end{picture}

\noindent Therefore, in view of \eqref{RE}, the number
$|\E_R(2)|$ of 2-cells in $\E_R$ can be estimated by
\begin{align*} |\E_R(2)| \le 1 +
\sum_{i=0}^{j-1} 3\cdot 2^{i+3} <  3\cdot 2^3 \cdot 2^j \ .
\end{align*}

Let $\D'$ be  a van Kampen  diagram over the presentation \eqref{pr1}
such that $\ph(\p \D') \equiv U$ and the $h$-tuple $\tau' = (\tau'_0,
\tau'_1, \dots)$, defined for $U$ by means of $\D'$, satisfies the
inequalities \eqref{tau1}, \eqref{tau2} in which $x = |U|$. As we saw in
the proof of Theorem~1.8, any word  $U \in \N(\infty)$ has such an
$h$-tuple $\tau'$, see \eqref{tau1}, \eqref{tau2}. Making use of diagrams
$\E_R$, we can turn $\D'$ into a diagram $\D$ over the presentation
\eqref{pr1a} of $\Ga_t$ so that $\ph(\p \D) \equiv U$ and
\begin{align*}
 |\D(2)| & \le  \tau'_0 + \sum_{j=1}^{[\log_2x]} 3\cdot 2^3 \cdot
 2^j\cdot \tau'_j \\ & \le  4C x^2 \log_2x +3\cdot 2^3\sum_{j=1}^{[\log_2x]}
 \frac{2^j}{2^{j-1}}C x^2 < 2^6 C x^2\log_2x \ .
\end{align*}
Hence, there is a diagram  $\D_W$  over \eqref{pr1a} such that
$\ph(\p \D_W) \equiv W$  and, in view of inequalities
\eqref{dU}, \eqref{dU2},
\begin{align*}
 |\D_W(2)| & \le 12 \cdot 2^{|W|/2} +  2^6 C |U|^2\log_2|U|  \\
           & \le 12\cdot 2^{|W|/2} + 2^6 \cdot 6^2 C \cdot  2^{|W|}\cdot  \log_2(6 \cdot 2^{|W|/2}   )
           \\  & \le C'_{2} |W|\cdot 2^{|W|}  \ ,
\end{align*}
where $C'_{2} >1$ is a constant. Since the length of every
relator of the presentation \eqref{pr1a} does not exceed 24, it
follows that
$$
|\D_W(j)| \le 24 |\D_W(2)| + \max(|W|, 1) \le \max (25C'_{2}
|W|\cdot 2^{|W|}, \ 1)
$$
for every $j =0,1,2$.   This proves that $f_j(x) \le 25C'_{2} x
 2^{x} = C_2 x  2^{x}$, as required. \qed

\section{Proofs of Theorem~1.10 and Corollary~1.11}

{\em Proof of Theorem~1.10:} Recall that $n \ge 2^{48}$ is a fixed
integer, $n$ is either odd or divisible by $2^9$ and $m \ge 2$. Under
these assumptions, lemmas of article \cite{Iv94} apply to diagrams over
the presentation $B(m,n, \infty) = \langle \A \, \| \, \R_B  \rangle$,
see \eqref{pr2}, and yield that, for every $i \ge 1$, the word $A_i$
exists and the limit group, defined by $B(m,n, \infty)$, is naturally
isomorphic to the free $m$-generator Burnside group $B(m,n) =
F(\A)/F(\A)^n$, where $F(\A)$ is the free group over $\A$.

Let $W$ be a word in the normal closure $\langle \langle \R_B
\rangle\rangle$  of $\R_B$. Then there exists a reduced diagram
$\D$ over $B(m,n, \infty)$, and hence over $B(m,n, i_0)$ for
some $i_0 \ge 0$, see  \eqref{pr2i}, such that $\ph (\p \D)
\equiv W$. Recall that the boundary $\p \D$ of a  van Kampen
(or disk) diagram $\D$ is oriented clockwise and the boundary
$\p \Pi$ of a face of $\D$ is oriented counterclockwise. Here
and below we are using the notation and terminology of article
\cite{Iv94} and the reader is referred to  \cite{Iv94} for more
details.

By induction on the perimeter $|\p \D|$ of a reduced diagram
$\D$ over  $B(m,n, i_0)$, we will be proving that the number
$| \D(1) |$ of 1-cells  in $\D$ does not exceed $|\p \D|^\ff$.
Since this claim is trivial when $\D$ contains no 2-cells, in
which case  $| \D(1) | = |\p \D|/2$, we may assume that $\D$
contains a 2-cell.

Denote $w {}= \p \D$, $x {}= |\p \D|$  and consider a
factorization $w = w_1 w_2 \dots w_8$, where for every $i=1,
\dots, 8$
\begin{equation}\label{w}
x/8 -1 = |w|/8 -1 < |w_i| <  |w|/8 +1 = x/8 +1
 \ ,
\end{equation}
where $|p|$ denotes the length of a path $p$.

 By Lemmas 5.7, 9.8 \cite{Iv94},  we can find a face
$\Pi$ and a system of subdiagrams $\Ga_1, \dots, \Ga_8$ in
$\D$, see Figure~4,  such that $\Ga_i$ is a contiguity
subdiagram between $\Pi$ and $w_i$ and
\begin{equation}\label{th}
\sum_{i=1}^8 | \Ga_i \wedge \p \Pi | > \theta |\p \Pi | \ ,
\quad \text{where} \quad \theta = 0.99 \  .
\end{equation}
We remark that, for some $i$, $\Ga_i$ might not be defined.
However, in view of  \eqref{th}, at least one of $\Ga_1, \dots,
\Ga_8$ is  defined.

Let $\p \Ga_i =  b_i u_i c_i q_i$ be the standard boundary of
$\Ga_i$ (if $\Ga_i$ is defined), where $u_i = \Ga_i \wedge
w_i$, $q_i = \Ga_i \wedge \p \Pi$,  and let $w_i = r_i u_i
s_i$, $i = 1, \dots, 8$, see Figure~4. If both $\Ga_i$ and
$\Ga_{i+1}$ are defined, here and below indices are considered
mod$\, 8$, then we let $q_{i+1} t_i q_i$ be a subpath of $\p
\Pi$ and $\D_i$ denote the subdiagram of $\D$ bounded by the
path $\p \D_i {}= c_i^{-1} s_i r_{i+1} b_{i+1}^{-1} t_i$, see
Figure~4. Informally, $\D_i$ sits between $\Ga_i$ and
$\Ga_{i+1}$ in the annulus $\D - \mbox{Int}\, \Pi$ and $t_i$
sits between $q_{i+1}$ and $q_{i}$ in the cycle $\p \Pi$.

\unitlength .6mm \linethickness{1pt}
\ifx\plotpoint\undefined\newsavebox{\plotpoint}\fi 
\begin{picture}(195.25,110)(8.00,-14.00)
\put(19.75,32.50){\line(1,0){166.75}}
\put(42.50,62.50){\line(0,-1){30.00}}
\put(82.50,62.50){\line(0,-1){30.00}}
\put(132.50,62.50){\line(0,-1){30.00}}
\put(172.50,62.50){\line(0,-1){30.00}}
\put(24.00,62.50){\vector(1,0){6.50}}
\put(18.50,62.50){\line(1,0){168.00}}
\put(48.75,62.50){\vector(1,0){3.75}}
\put(69.50,62.50){\vector(1,0){5.00}}
\put(102.00,62.50){\vector(1,0){7.50}}
\put(139.50,62.50){\vector(1,0){4.00}}
\put(160.00,62.50){\vector(1,0){5.50}}
\put(180.50,62.50){\vector(1,0){3.50}}
\put(42.50,51.00){\vector(0,-1){4.25}}
\put(82.50,45.75){\vector(0,1){3.50}}
\put(132.50,49.75){\vector(0,-1){4.00}}
\put(172.50,45.50){\vector(0,1){4.25}}
\put(185.25,32.50){\vector(-1,0){3.75}}
\put(159.50,32.50){\vector(-1,0){6.50}}
\put(110.75,32.50){\vector(-1,0){5.25}}
\put(65.50,32.50){\vector(-1,0){4.50}}
\put(31.75,32.50){\vector(-1,0){3.75}}
\put(61.75,62.50){\circle*{1.5}}
\put(152.75,62.50){\circle*{1.5}}
\put(186.00,62.50){\line(1,0){8.75}}
\put(186.50,32.50){\line(1,0){8.75}} \put(106.50,65.50){$u_i$}
\put(47,65.50){$s_{i-1}$} \put(72,65.50){$r_i$}
\put(26,65.50){$u_{i-1}$} \put(44.50,48.50){$c_{i-1}$}
\put(27.75,27){$q_{i-1}$} \put(61.75,27){$t_{i-1}$}
\put(106.75,27){$q_{i}$} \put(154.50,27){$t_{i}$}
\put(85.00,47){$b_{i}$} \put(134.50,47.25){$c_{i}$}
\put(174.50,46.75){$b_{i+1}$} \put(141.75,65.50){$s_{i}$}
\put(158.00,65.50){$r_{i+1}$} \put(180,65.50){$u_{i+1}$}
\put(182.00,27){$q_{i+1}$} \put(107.00,47.25){$\Gamma_i$}
\put(61.00,48.00){$\D_{i-1}$} \put(25.00,49.00){$\Gamma_{i-1}$}
\put(151.75,48.25){$\D_i$} \put(188.50,48.50){$\Gamma_{i+1}$}
\put(81.75,75.25){$w_i = r_i u_i s_i$}
\put(122.00,17.75){$\Pi$}
\put(128.25,80.00){\vector(1,0){17.25}}
\put(62.75,18.25){\vector(-1,0){19.50}} \put(130.50,83.50){$\p
\D = \dots  r_i u_i s_i  \dots$} \put(38,10.00){$\p \Pi = \dots
t_i q_i t_{i-1}q_{i-1}  \dots $} \put(110.00,0){Figure 4}
\end{picture}

It follows from Lemmas~3.1, 9.8 \cite{Iv94} applied to $\Ga_i$
that
\begin{equation}\label{ga}
\max(|b_i|, |c_i|) < \gamma | \p \Pi |  \ , \quad \text{where}
\quad \gamma = 2^{-33}  \ , \quad i =1, \dots, 8 \ .
\end{equation}

By Lemmas~6.1, 9.8 \cite{Iv94} applied to $\Ga_i$, we also have
that
\begin{equation}\label{ro1}
\rho |q_i | \le |u_i| + |b_i|+ |c_i| \ , \quad \text{where}
\quad \rho = 0.95  \ , \quad i =1, \dots, 8 \ .
\end{equation}

It follows from Lemmas~6.2, 9.8 \cite{Iv94} applied to the
diagram $\D$ and its face $\Pi$  that
\begin{equation}\label{ro2}
\rho n  \le \rho |\p \Pi | \le |\p \D | = x \ .
\end{equation}

Similar to Storozhev's arguments \cite[Section 28.2]{Ol89}, we
consider two cases which are whether or not all of the
subdiagrams  $\Ga_1, \dots, \Ga_8$ are actually defined.
\medskip

First we assume that at least one of $\Ga_1, \dots, \Ga_8$ is
missing. Reindexing if necessary, suppose that $\Ga_1, \Ga_{j}$
are present and $\Ga_2, \dots, \Ga_{j-1}$ are missing, where
$j-1 \in \{ 3, \dots, 8\}$. Let us emphasize that it { is}
possible that  $j-1 =8$ and $j=1$, that is, the only contiguity
subdiagram, which is defined, is $\Ga_1$.

If $t$ is a path in $\D$, then $t_-$, $t_+$ denote the initial,
terminal, respectively, vertices of $t$. Consider a path $p =
c_1 v^{-1} b_j$ in $\D$, where $v$ is the subpath of $\p
\Pi^{-1}$ that connects the vertices $(c_1)_+$, $(b_j)_-$  and
contains no paths $q_i$, $i =1, \dots, 8$. By estimates
 \eqref{ga},
\begin{equation}\label{p}
|p|  < (1-\theta + 2\gamma)  | \p \Pi |  \ .
\end{equation}

The vertices $p_-, p_+ \in \p \D$ define a factorization $\p \D
= y_1 y_2$, where $y_2$ contains the subpath $u_1$ of $\p \D$.
Cutting $\D$ along $p$, we obtain two diagrams $\E_1$, $\E_2$
with $\p \E_1 = y_1 p^{-1}$, $\p \E_2 = p y_2$.  Since $y_1$
contains $w_2$, it follows from \eqref{w} and \eqref{ro2} that
\begin{equation}\label{y1}
|y_1|  \ge x/8 -1 \ge  (1/8 - (\rho n)^{-1}) x \ge (\rho/8
-n^{-1}) | \p \Pi | \ .
\end{equation}

Adding up all estimates \eqref{ro1} (for those  $\Ga_i$ that
are defined), in view of \eqref{th} and \eqref{ga}, we have
\begin{equation*}
 \rho\theta |\p \Pi |  <  \sum_{i=1}^8 |u_i| + 14\gamma |\p  \Pi |
 \le |y_2| + 14\gamma |\p  \Pi |   \ ,
\end{equation*}
because every $u_i = \Gamma_i \wedge w_i = \Gamma_i \wedge \p
\D$ is a subpath of $y_2$. Therefore,
\begin{equation}\label{y2}
 |y_2|  >   (\rho\theta - 14\gamma) |\p \Pi | \ .
\end{equation}

Now we see from  \eqref{p}, \eqref{y1}, \eqref{y2} that
\begin{align}\notag
|p|  & <  (1-\theta + 2\gamma) \min\{ (\rho/8 -n^{-1})^{-1}
|y_1|, \
 (\rho\theta - 14\gamma)^{-1}  |y_2|\} \\ & < 0.1\min(|y_1|, |y_2|) \ . \label{011}
\end{align}
In particular, $|\p \E_1 |, |\p \E_2 | < x$. Since $ \E_1 ,
\E_2$ are reduced diagrams,   the induction hypothesis applies
to $\E_1 , \E_2$  and yields that $| \E_k(1) | \le |\p
\E_k|^\ff$, $k=1,2$. Since $|\D(1)| \le |\E_1(1)| + |\E_2(1)|$,
it follows that
\begin{align}\notag
|  \D(1) |   & \le     ( |p| + |y_1| )^\ff +  ( |p| + |y_2|
)^\ff
\\   \notag  & \le \max(|p| + |y_1|, |p| + |y_2| )^{7/12}( |y_1|+ |y_2|+ 2|p|
) \\ & \le \max(|p| + |y_1|, |p| + |y_2| )^{7/12}( x+ 2|p|) \ .
\label{d1}
\end{align}
In view of estimate  \eqref{011},
\begin{align}\notag
 \max(|p| + |y_1|, |p| + |y_2| ) & = \max(x - |y_2| + |p| , x - |y_1| + |p|  )
 \\  & \le x - 10|p| + |p| =   x - 9|p|   \ .
\label{d2}
\end{align}

By inequalities  \eqref{ro2} and \eqref{p},
\begin{equation}\label{d3}
|p|  < \rho^{-1}(1-\theta + 2\gamma)x < 0.02x   \ .
\end{equation}

Let $r = \frac{|p|}{x}$. Then $0 < r < 0.02$ and we can derive
from \eqref{d1}, \eqref{d2}, \eqref{d3} that
\begin{align*}
|  \D(1) | & < (1-9r)^{7/12}(1+2r) x^{19/12} <
(1-9r)^{1/2}(1+2r) x^{19/12}\\ & \le ((1-9r)(1+2r)^2)^{1/2}
x^{19/12} < ((1-9r)(1+6r))^{1/2} x^{19/12} < x^{19/12} \ ,
\end{align*}
as required.
\medskip

Now suppose that all of $\Ga_i$, $i =1, \dots, 8$, { are}
defined. Taking $\text{Int}\, \Pi$ out of $\D$, we get an
annular diagram $\D_0$. Cutting  $\D_0$ along all paths $b_i,
c_i$, we obtain 16 disk  diagrams $\Ga_i, \D_i$, $i =1, \dots,
8$, see Figure~4. It follows from the definitions that
\begin{align*}
| \p  \Ga_i |  = |b_i| + |u_i| + |c_i| + |q_i| \ , \  \quad |
\p \D_i |  \le |w_i| + |w_{i+1}| +  |c_i|+ |b_{i+1}| + |t_i| \
.
\end{align*}

Hence, in view of estimates \eqref{w}, \eqref{th}, \eqref{ga},
\eqref{ro1}, \eqref{ro2}, we further have
\begin{align} \notag
| \p  \Ga_i | & < (1+ \rho^{-1})( |u_i| + 2 \gamma   |\p \Pi |)
\\ &  \le  (1+ \rho^{-1})( 1/8 +(\rho n)^{-1} +
2 \gamma  \rho^{-1}) x  < 0.258 x  \ , \label{Gi}
\\ | \p \D_i |  &  \le 2 (1/8 +(\rho n)^{-1} )x +
(1 - \theta + 2 \gamma )|\p \Pi
|  \notag \\  &  \le (1/4 +2(\rho n)^{-1} + \rho^{-1}(1 -
\theta + 2 \gamma) )x < 0.262 x \ .\label{Di}
\end{align}

Since $\Ga_i$, $\D_i$, $i =1, \dots, 8$, are reduced diagrams,
the induction hypothesis applies to them and yields that $|
\Ga_i(1) | \le | \p \Ga_i |^\ff$, $| \D_i(1) | \le | \p \D_i
|^\ff$, $i =1, \dots, 8$. Therefore,
\begin{align} \notag
| \D(1) | & \le \sum_{i=1}^8  (  | \Ga_i(1) | + | \D_i(1) | )
\le \sum_{i=1}^8  (  | \p \Ga_i |^\ff + | \p \D_i |^\ff ) \\
& \le  \max( | \p \Ga_i |, | \p \D_i | )^{7/12} \cdot
\sum_{i=1}^8 (  | \p \Ga_i | + | \p \D_i | ) \ . \label{sum}
\end{align}
It follows from the definitions and estimates \eqref{ga},
\eqref{ro2} that
\begin{align} \notag
\sum_{i=1}^8 (  | \p \Ga_i | + | \p \D_i | ) & = |w| + |\p \Pi
|  +   \sum_{i=1}^8 (|b_i| + |c_i|)  \le   ( 1 +\rho^{-1})x +
32\gamma  |\p  \Pi |    \\ & \le   (1 +\rho^{-1} + 32\gamma
\rho^{-1}  ) x < 2.06 x \ . \label{sum2}
\end{align}

In view of estimates \eqref{Gi}, \eqref{Di}, \eqref{sum2}, it
follows from \eqref{sum} that
\begin{align*}
| \D(1) |  \le (0.262x)^{7/12}\cdot 2.06 x < x^\ff \ ,
\end{align*}
as desired. Thus the desired upper bounds  $f_{1,B}(x) \le
x^\ff $, $f_{1,B(i_0)}(x) \le x^\ff $ are proven. The
inequalities  $f_{0,B}(x) \le 2 x^\ff$, $f_{0,B(i_0)}(x) \le 2
x^\ff$  now follow from Theorem~1.4(a) and the inequalities
 $f_{2,B}(x) \le \tfrac 2n x^\ff $,  $f_{2,B(i_0)}(x) \le \tfrac
2n x^\ff $ hold because $|R| \ge n$ for every relator $R$ of $B(m,n,
\infty)$, $B(m,n,i_0)$. The proof of Theorem~1.10 is complete. \qed

\medskip

{\em Proof of Corollary~1.11:} Let $W$ be a word over the alphabet
$\A^{\pm 1}$. Then $W$ represents the identity element of the group
$B(m,n) = F(\A)/F(\A)^n$, defined by $B(m,n, \infty)$, if and only if
there exists a van Kampen diagram $\D_W$ over $B(m,n, \infty)$ with
$\ph(\p \D_W) \equiv W$. By Theorem~1.10, we can assume that $|\D_W(1)|
\le |W|^\ff$. As in the proof of Theorem~1.3, this polynomial  bound
enables us to take a van Kampen diagram $\D$ with $|\D(1)| \le |W|^\ff$
as a certificate to verify whether $W = 1$ in $B(m,n)$. In polynomial
time with respect to $|W|$, we can certify that $\ph(\p \D) \equiv W$
and, for every face $\Pi$ in $\D$, the label $\ph(\p \Pi)$ is the $n$th
power of a word. If these conditions are satisfied, then $W =1$ in
$B(m,n)$, as required.

Now let $U$, $V$ be two words over the alphabet $\A^{\pm 1}$. It follows
from Lemmas 6.3, 9.2 \cite{Iv94} and Theorem~1.10 that $U$, $V$ are
conjugate in $B(m,n)$ if and only if there exists an annular diagram
$\D_{U,V}$ over $B(m,n, \infty)$ such that the oriented components $p$,
$q$ of the boundary $\p \D_{U,V}$ are labelled by the cyclic words $U$,
$V^{-1}$ and
$$
|\D_{U,V}(1)| \le (1.01(|U|+|V|))^\ff \le 1.02(|U|+|V|)^\ff  \
.
$$
Hence, as above, an annular diagram $\D$ with $|\D(1)| \le
1.02(|U|+|V|)^\ff$ can be used  to certify, in polynomial time
of $|U|+|V|$,  that  $U$ and $V$ are conjugate in $B(m,n)$ by
checking that $\ph(p) \equiv U$, $\ph(q) \equiv V^{-1}$ and
that, for every face $\Pi$ in $\D$, $\ph(\p \Pi)$ is the $n$th
power of a word. \qed

\medskip

{\em Acknowledgements:}  We wish to thank  M.~Sapir for helpful
discussions and for pointing out to the notion of the verbal Dehn
function  and to  other results of \cite{OS2}. We also thank C.~Jockush
and I.~Kapovich for providing Example~2.3. We are indebted to P. de la
Harpe for many useful suggestions and to the  anonymous referee for a
number of valuable remarks and for providing Example~2.4.


\begin{thebibliography}{[75]}

\bibitem[1]{A75}
S. I. Adian,  {\em The Burnside problem and identities in
groups},  Nauka, Moscow, 1975; English translation:
Springer-Verlag, 1979.



\bibitem[2]{BMS}
G. Baumslag, C. F. Miller III, and H. Short, {\em Isoperimetric
inequalities and the homology of groups}, {Invent. Math.} {\bf
113}(1993), 531--560.

\bibitem[3]{Bir} J.-C. Birget, {\em  Infinite string rewrite systems and complexity},
J. Symbolic Comp. {\bf 25}(1998) 759--793.

\bibitem[4]{Br} S. G. Brick, {\em On Dehn functions and products of groups},
Trans. Amer. Math. Soc. {\bf 335}(1993), 369--384.


\bibitem[5]{BORS}
 J.-C. Birget,   A. Yu. Ol'shanskii, E. Rips, and M. V. Sapir,
 {\em Isoperimetric functions of groups and computational
complexity of the word problem}, {Ann. of Math.} {\bf
156}(2002), 467--518.


\bibitem[6]{BB} M. R. Bridson and N. Brady, {\em  There is only one gap in the
isoperimetric spectrum}, Geom. Funct. Anal. {\bf 10}(2000)
1053--1070.


\bibitem[7]{G80}
R. I. Grigorchuk,  {\em On  the Burnside problem for periodic
groups}, {Funct. Anal. Appl.} {\bf 14}(1980), 53--54.


\bibitem[8]{G84}
R. I. Grigorchuk,  {\em The growth degrees of finitely
generated groups and the theory of invariant means}, {Izv.
Akad. Nauk. SSSR. Ser. Mat.} {\bf 48}(1984), 939--985.



\bibitem[9]{G98}
R. I. Grigorchuk, {\em An example of a finitely presented
amenable group that does not belong to the class EG}, {Mat.
Sbornik} {\bf 189}(1998), 79--100.


\bibitem[10]{G99}
 R. I. Grigorchuk, {\em On the system of defining relations and the
Schur multiplier of periodic groups generated by finite
automata}, in  {\em  Groups St. Andrews 1997 in Bath, I},
London Math. Soc. Lecture Note Ser. {\bf 260}, 1999,  290--317.


\bibitem[11]{G00}
 R. I. Grigorchuk, {\em Just infinite branch groups}, in
 {\em New horizons in pro-$p$-groups},  Progr. Math. {\bf 184}, Birkh\"auser,
2000, 121--179.



\bibitem[12]{G05}
 R. I. Grigorchuk,  {\em  Solved and unsolved problems around one
group}, in {\em Infinite groups: geometric, combinatorial and
dynamical aspects}, Progr. Math. {\bf 248}, Birkh\"auser, 2005,
117--218.



\bibitem[13]{GH}
R. I. Grigorchuk and  P. de la Harpe, {\em  Limit behaviour of
exponential growth rates for finitely generated groups}, In
{\em Essays on geometry and related topics},  Monogr. Enseign.
Math. {\bf  38}, Enseignement Math., 2001, 351--370.




\bibitem[14]{Gromov}
M. Gromov,  {\em Hyperbolic groups}, In {\em Essays in Group
Theory} (S. Gersten, ed.), MSRI Publ. 8, Springer-Verlag, 1987,
75--263.



\bibitem[15]{Guba}
V. S. Guba, {\em The Dehn function of Richard Thompson's group
$F$ is quadratic}, Invent. Math.  {\bf 163}(2006), 313--342.


\bibitem[16]{GS}
V. S. Guba and M. V. Sapir, {\em On Dehn functions of free
products of groups}, {Proc.  Amer. Math. Soc. } {\bf
127}(1999), 1885--1891.

\bibitem[17]{harpe}
P. de la Harpe,  {\em Topics in geometric group theory},
 Univ. Chicago Press, 2000.


\bibitem[18]{Iv92}
S. V. Ivanov, {\em On the Burnside problem on periodic groups},
{Bull. Amer. Math. Soc.} {\bf 27}(1992), 257--260.

\bibitem[19]{Iv94}  S. V. Ivanov,  {\em The free Burnside groups of sufficiently
large exponents}, {Internat. J. Algebra Comp.}  {\bf
4}\,(1994), 1--308.

\bibitem[20]{Iv98}
 S. V. Ivanov, {\em On the Burnside problem for groups of even exponent},
{Documenta Mathematica, Extra Volume ICM-98} {\bf II}(1998),
67--76.

\bibitem[21]{Iv05} S. V. Ivanov,  {Embedding free Burnside groups in
finitely presented groups},  {\em Geom. Dedicata} {\bf
111}(2005), 87--105.

\bibitem[22]{Iv06}  S. V. Ivanov,
{\em On balanced presentations of the trivial group}, {Invent.
Math.}  {\bf 165}(2006), 525--549.


\bibitem[23]{LS} R. C. Lyndon and P. E. Schupp, {\em Combinatorial
group theory}, Springer-Verlag,  1977.


\bibitem[24]{Lys85}
 I. G. Lysenok,  {\em A set of defining relations for the Grigorchuk group},
{Matem. Zametki} {\bf 38}(1985), 503--516.


\bibitem[25]{Lys96}
 I. G. Lysenok,   {\em Infinite Burnside groups of even period},
{Izv. Ross. Akad. Nauk Ser. Mat.} {\bf 60}(1996), 3--224.


\bibitem[26]{MaOt} K. Madlener, F. Otto,  {\em Pseudo-natural algorithms
for the word problem for finitely presented  groups}, { J.
Symbolic Comp.} {\bf 1}(1985), 383--418.


\bibitem[27]{MKS}   W. Magnus,  J.  Karras,  D. Solitar,
{\em Combinatorial group theory}, Interscience Publ., 1966.



\bibitem[28]{NA}
P. S. Novikov and S. I. Adian, {\em  On infinite periodic
groups I, II, III}, Izv. Akad. Nauk SSSR Ser. Mat. {\bf
32}(1968), 212--244, 251--524, 709--731.

\bibitem[29]{Ol82}
A. Yu. Ol'shanskii, {\em On the Novikov-Adian theorem}, {Mat.
Sbornik} {\bf 118}(1982), 203--235.

\bibitem[30]{Ol89}
A. Yu. Ol'shanskii,  {\em Geometry of defining relations in
groups}, Nauka, Moscow, 1989; English translation: {\em Math.
and Its Applications, Soviet series},  vol. 70, Kluwer Acad.
Publ., 1991.

\bibitem[31]{Ol91}
A. Yu. Ol'shanskii, {\em   Hyperbolicity of groups with
subquadratic isoperimetric inequality}, {Internat. J. Algebra
Comput.} {\bf 1}(1991), 281--289.

\bibitem[32]{OS2}  A. Yu. Ol'shanskii and M. V. Sapir,
{\em Embeddings of relatively free groups into finitely
presented groups}, Contemp. Math. {\bf   264}(2000), 23--47.


\bibitem[33]{OS1} A. Yu. Ol'shanskii and M. V. Sapir,
{\em Length and area functions on groups and quasi-isometric
Higman embeddings}, Internat. J. Algebra  Comput. {\bf
11}(2001), 137--170.


\bibitem[34]{OS} A. Yu. Ol'shanskii and M. V. Sapir,
{\em Non-amenable finitely presented torsion-by-cyclic groups},
Publ. Math. IHES {\bf  96}(2002), 43--169.


\bibitem[35]{OS4} A. Yu. Ol'shanskii and M. V. Sapir,
{\em Groups with small Dehn functions and bipartite chord
diagrams}, Geom. Funct. Anal. {\bf 16}(2006), 1324--1376.

\bibitem[36]{Pap}
C. H.~Papadimitriou, \emph{Computational complexity},
Addison-Wesley Publ., 1994.

\end{thebibliography}
\end{document}